\documentclass[final]{siamltex}
\usepackage{amsmath,amsfonts,graphicx,color}
\newtheorem{remark}[theorem]{Remark}

\newlength\symbolwidth

\newenvironment{symbollist}{\list{}{\labelwidth\symbolwidth%
    \leftmargin\labelwidth\advance\leftmargin by\labelsep%
     \topsep %
     \parsep %
    \itemsep \parsep%
}}{\endlist}

\newcommand{\real}{{\mathbb{R}}}

\renewcommand{\natural}{{\mathbb{N}}}
\newcommand{\eps}{\epsilon}
\newcommand{\ov}{\overline}

\newcommand{\pder}[2]{\frac{\partial #1}{\partial #2}}

\newcommand{\setdef}[2]{\left\{#1 \; | \; #2\right\}}
\newcommand{\map}[3]{#1: #2 \rightarrow #3}
\newcommand{\versus}[1]{\operatorname{vrs}(#1)}
\newcommand{\until}[1]{\{1,\dots,#1\}}

\newcommand{\setLieder}[2]{\widetilde{\Lc}_{#1} #2}
\newcommand{\proj}[1]{\operatorname{proj}_{#1}}

\newcommand{\co}{\operatorname{co}}
\newcommand{\interior}[1]{\operatorname{int}(#1)}
\newcommand{\sm}{\operatorname{sm}} 
\renewcommand{\lg}{\operatorname{lg}}
\newcommand{\CC}{\operatorname{CC}} 
\newcommand{\IC}{\operatorname{IC}}
\newcommand{\CR}{\operatorname{CR}} 
\newcommand{\IR}{\operatorname{IR}}
 
\newcommand{\D}{\operatorname{D}}
\newcommand{\DC}{\operatorname{DC}} 
\newcommand{\SP}{\operatorname{SP}}
\newcommand{\ED}{\operatorname{Ed}} 
\newcommand{\VE}{\operatorname{Ve}}
\newcommand{\LN}{\operatorname{Ln}} 
\newcommand{\Id}{\operatorname{Id}}

\newcommand{\Lc}{{\mathcal{L}}}
\newcommand{\HHDC}{{\mathcal{H}}_{\text{DC}}}
\newcommand{\HHSP}{{\mathcal{H}}_{\text{SP}}}

\newcommand{\VV}{{\mathcal{V}}}
\newcommand{\WW}{{\mathcal{W}}}
\newcommand{\NN}{{\mathcal{N}}}
\newcommand{\GG}{{\mathcal{G}}}
\newcommand{\FF}{{\mathcal{F}}}

\renewcommand{\theenumi}{(\roman{enumi})}

\title{Coordination and geometric optimization \\ via distributed dynamical
  systems\footnotemark[1]}

\author{Jorge Cort\'es\footnotemark[2] \and Francesco
  Bullo\footnotemark[2]}

\begin{document}
\maketitle \renewcommand{\thefootnote}{\fnsymbol{footnote}}

\footnotetext[1]{Submitted to the SIAM Journal on Control and Optimization
  on May 27, 2003.
  %% This version: \today. 
  This work was supported by DARPA/AFOSR MURI Award F49620-02-1-0325. A
  preliminary version of this manuscript has been submitted to the 2003
  IEEE Control and Decision Conference, Maui, Hawaii.}

\footnotetext[2]{Coordinated Science Laboratory, University of Illinois at
  Urbana-Champaign, 1308 W. Main St., Urbana, IL 61801, United States, Ph.  +1
  217 244-8734 and +1 217 333-0656, Fax. +1 217 244-1653,
  \texttt{\{jcortes,bullo\}@uiuc.edu},
  \texttt{http://motion.csl.uiuc.edu/\~{}\{jorge,bullo\}}}
\renewcommand{\thefootnote}{\arabic{footnote}}

\begin{abstract}
  This paper discusses dynamical systems for disk-covering and
  sphere-packing problems.  We present facility location functions from
  geometric optimization and characterize their differentiable properties.
  We design and analyze a collection of distributed control laws that are
  related to nonsmooth gradient systems. The resulting dynamical systems
  promise to be of use in coordination problems for networked robots; in
  this setting the distributed control laws correspond to local
  interactions between the robots.  The technical approach relies on
  concepts from computational geometry, nonsmooth analysis, and the
  dynamical system approach to algorithms.
\end{abstract}


\pagestyle{myheadings} \thispagestyle{plain}
 
\markboth{Jorge Cort\'es and Francesco Bullo}%
{Coordination and geometric optimization via distributed dynamical systems}
%% shorter: {Distributed coordination and geometric optimization}

\section{Introduction}\label{sec:introduction}
Consider $n$ sites $(p_1,\dots,p_n)$ evolving within a convex polygon $Q$
according to one of the following interaction laws: (i) each site moves away
from the closest other site or polygon boundary, (ii) each site moves toward
the furthest vertex of its own Voronoi polygon, or (iii) each site moves
toward a geometric center (circumcenter, incenter, centroid, etc) of its own
Voronoi polygon.  Recall that the Voronoi polygon of the $i$th site is the
closed set of points $q\in Q$ closer to $p_i$ than to any other $p_j$.

These and related interaction laws give rise to strikingly simple dynamical
systems whose behavior remains largely unknown. What are the critical
points of such dynamical systems? What is their asymptotic behavior? Are
these systems optimizing any aggregate function? In what way do these local
interactions give rise to distributed systems?  Does any biological
ensemble evolve according to these behaviors and are they of any
engineering use in coordination problems?  These are the questions that
motivate this paper.

% are these dynamical systems gradient algorithms for some functions in
% some sense? 

\subsection*{Coordination in robotics, control, and biology}
Coordination problems are becoming increasingly important in numerous
engineering disciplines.  The deployment of large groups of autonomous
vehicles is rapidly becoming possible because of technological advances in
computing, networking, and miniaturization of electro-mechanical systems.
These future multi-vehicle networks will coordinate their actions to
perform challenging spatially-distributed tasks (e.g., search and recovery
operations, exploration, surveillance, and environmental monitoring for
pollution detection and estimation).  This future scenario motivates the
study of algorithms for autonomy, adaptation, and coordination of
multi-vehicle networks.  It is also important to take into careful
consideration all constraints on the behavior of the multi-vehicle network.
Coordination algorithms need to be adaptive and distributed in order for
the resulting closed-loop network to be scalable, to comply with bandwidth
limitations, to tolerate failures, and to adapt to changing environments,
topologies and sensing tasks. The interaction laws introduced above have
these properties and, remarkably, they optimize network-wide performance
measures for meaningful spatially-distributed tasks.

Coordinated group motions are also a widespread phenomenon in biological
systems.  Some species of fish spend their lives in schools as a defense
mechanism against predators. Others travel as swarms in order to protect an
area that they have claimed as their own. Flocks of birds are able to
travel in large groups and act as one unit. Other animals exhibit
remarkable collective behaviors when foraging and selecting food.  Certain
foraging behaviors include individual animals partitioning their
environment in nonoverlapping individual zones whereas other species
develop overlapping team areas.  These biological network systems possess
extraordinary dynamic capabilities without apparently following a group
leader.  Yet, these complex coordinated behaviors emerge while each
individual has no global knowledge of the network state and can only plan
its motion according to the observation of its closest neighbors.

\subsection*{Facility location,  nonsmooth stability analysis and cooperative
  control} 

To analyze the interaction laws introduced above we rely on concepts and
methods from various disciplines.  Facility location problems play a
prominent role in the field of geometric
optimization~\cite{PKA-MS:98,VB-HM-VS:99}. Facility location pervades a
broad spectrum of scientific and technological areas, including resource
allocation (where to place mailboxes in a city or cache servers on the
internet), quantization and information theory, mesh and grid optimization
methods, clustering analysis, data compression, and statistical pattern
recognition.  Smooth multi-center functions for so-called centroidal
Voronoi configurations and smooth distributed dynamical systems are
presented in~\cite{JC-SM-TK-FB:02j,QD-VF-MG:99}.  Multi-center functions
are studied in resource allocation problems~\cite{ZD:95,AS-ZD:96} and in
quantization theory~\cite{RMG-DLN:98,SPL:82}.  The role of Voronoi
tessellations and computational geometry in facility location is discussed
in~\cite{AO-BB-KS-SNC:00,JMR-GTT:90}.

The notion and computational properties of the generalized gradient are
throughly studied in nonsmooth analysis~\cite{FHC:83}. In particular, tools
for establishing stability and convergence properties of nonsmooth
dynamical systems are presented in~\cite{AB-FC:99,AFF:88,DS-BP:94}.
Finally, we refer to~\cite{RWB:91,UH-JBM:94} for guidelines on how to
design dynamical systems for optimization purposes, and
to~\cite{DPB-JNT:97} for gradient descent flows in distributed computation
in settings with fixed-communication topologies.

Recent years have witnessed a large research effort focused on motion
planning and formation control problems for multi-vehicle
systems~\cite{JPD-JPO-VK:01,AJ-JL-ASM:02,NEL-EF:01,YL-KMP-MMP:00,ROS-RMM:03b,IS-MY:99,HT-AJ-GJP:03b}.
Within the literature on behavior-based robotics, heuristic approaches to
the design of interaction rules and emerging behaviors have been
investigated (see~\cite{RCA:98} and references therein).  Along this
specific line of research, no formal results guaranteeing the correctness
of the proposed algorithms or their optimality with respect to an aggregate
objective are currently available. The aim of this work is to design
distributed coordination algorithms for dynamic networks as well as to
provide formal verifications of their asymptotic correctness.  A key aspect
of our treatment is the inherent complexity of studying networks whose
communication topology changes along the system evolution, as opposed to
networks with fixed communication topologies. This key aspect is present in
the analysis of distributed control laws
in~\cite{AJ-JL-ASM:02,IS-MY:99,HT-AJ-GJP:03b} and of agreement protocols
in~\cite{ROS-RMM:03b}.

%% This work contributes to an ongoing research program initiated
%% in~\cite{JC-SM-TK-FB:02j} on cooperative control problems for multi-vehicle
%% networks

\subsection*{Statement of contributions}
We consider two facility location functions from geometric optimization
that characterize coverage performance criteria.  A collection of sites
provides optimal service to a domain of interest if (i) it minimizes the
largest distance from any point in the domain to one of the sites, or (ii)
it maximizes the minimum distance between any two sites.  In other words,
if $P=(p_1,\dots,p_n)$ are $n$ sites evolving within a convex polygon $Q$,
we extremize the \emph{multi-center functions}
\begin{align*}
  \max_{q\in Q} \left\{ \min_{i\in\{1,\ldots,n\}} d(q,p_i) \right\} \, ,
  \quad \min_{ i\not=j\in\until{n} } \left\{ {\textstyle \frac{1}{2}}
    d(p_i,p_j), d(p_i,\partial Q)\right\} \, ,
\end{align*}
where $d(p,q)$ and $d(p,\partial Q)$ are the distances between $p$ and $q$,
and between $p$ and the boundary of $Q$, respectively.  (The role of the
$\frac{1}{2}$ factor will become clear later.) We study the differentiable
properties of these functions via nonsmooth analysis.  We show the
functions are globally Lipschitz and regular, we compute their generalized
gradients, and we characterize their critical points.  Under certain
technical conditions, we show that the local minima of the first
multi-center function are so-called circumcenter Voronoi configurations,
and that these critical points correspond to the solutions of disk-covering
problems.  Similarly, under analogous technical conditions, we show that
the local maxima of the second multi-center function are so-called incenter
Voronoi configurations, and that these critical points correspond to the
solutions of sphere-packing problems.

Next, we aim to design distributed algorithms that extremize the multi-center
functions. Roughly speaking, by distributed we mean that the evolution of
each site depends at most on the location of its own Voronoi neighbors. We
study the generalized gradient flows induced by the multi-center functions
using nonsmooth stability analysis.  Although these dynamical systems possess
some convergence properties, they are not amenable to distributed
implementations. Next, drawing connections with quantization theory, we
consider two dynamical systems associated to each multi-center function.
First, we consider a novel strategy based on the generalized gradient of the
1-center functions of each site, and, second, we consider a geometric
centering strategy similar to the well-known Lloyd
algorithm~\cite{RMG-DLN:98,SPL:82}.

Remarkably, these strategies arising from the nonsmooth gradient
information have natural geometric interpretations and are indeed the local
interaction rule described earlier.  For the first (respectively second)
multi-center function, the first strategy corresponds to the interaction
law ``move toward the furthest vertex of own Voronoi polygon''
(respectively, ``move away from the closest other site or polygon
boundary'', and the second strategy corresponds to the interaction law
``move toward circumcenter of own Voronoi polygon'' (respectively ``move
toward incenter of own Voronoi polygon'').  We prove the uniqueness of the
solutions of the resulting distributed dynamical systems and we analyze
their asymptotic behavior using nonsmooth stability analysis, showing that
the active sites will approach the corresponding centers of their own
Voronoi cells.

% The gradient of the locational optimization function is a separable
% function and naturally gives rise to a variety of distributed
% flows~\cite{JC-SM-TK-FB:02j,QD-VF-MG:99,AO-BB-KS-SNC:00}.

Two of our results are related to well-known conjectures in the locational
optimization literature~\cite{ZD:95,AS-ZD:96}: (i) that the first
multi-center problem is equivalent to a disk covering problem (how to cover
a region with possibly overlapping disks of equal minimum radius), and (ii)
that the generalized Lloyd strategy ``move toward circumcenter of own
Voronoi polygon'' converges to the set of circumcenter Voronoi
configurations.

\subsection*{Organization}
The paper is organized as follows. Section~\ref{sec:preliminaries} provides
the preliminary concepts on Voronoi partitions, nonsmooth analysis, stability
analysis, and gradient flows, and introduces the multi-center problems.
Section~\ref{sec:1-center} presents a complete treatment on the functions
analysis and algorithm design for the 1-center problems.
Section~\ref{sec:multi-center-analysis} discusses the differentiable
properties and the critical points of the multi-center functions.
Section~\ref{sec:multi-center-design} introduces a number of dynamical
systems (smooth and nonsmooth, distributed and non-distributed) and analyzes
their asymptotic correctness.

% The resulting control laws correspond to basic
% interaction behaviors between the robots.

\section{Preliminaries and problem setup}\label{sec:preliminaries}

Let $\|\cdot\|$ denote the Euclidean distance function on $\real^N$ and let
$v \cdot w$ denote the scalar product of the vectors $v,w\in \real^N$. Let
$\versus{v}$ denote the unit vector in the direction of $0 \neq v \in
\real^N$, i.e., $\versus{v} = v / \| v \|$.  Given a set $S$ in $\real^N$, we
denote its convex hull by $\co(S)$ and its interior set by $\interior{S}$.
If $S$ is a convex set in $\real^N$, let $\map{\proj{S}}{\real^N}{S}$ denote
the orthogonal projection onto $S$ and let $\map{\D_S}{\real^N}{\real}$
denote the distance function to $S$.
% %
% \sindex{versus}{$\versus{v}$}{Unit vector in the direction of $0 \neq v \in \real^N$}%
% \sindex{proj}{$\proj{S}$}{Orthogonal projection onto the convex set $S$}%
% \sindex{distance}{$\D_S$}{Distance function to the convex set $S$}
% %
For $R>0$, let $\ov{B}_N(p,R) = \setdef{q \in \real^N}{\| p-q\| \le R}$,
and $B_N(p,R) = \interior{\ov{B}_N(p,R)}$.  A set $\{v_1,\dots,v_M\}$ of
vectors in $\real^N$ \emph{positively spans $\real^N$} if any $w \in
\real^N$ can be written as $w=\sum_{l=1}^M a_l v_l$, with $a_l \ge 0$, $l
\in \until{M}$.  The following simple lemma, e.g., see~\cite{HC:99},
characterizes this situation.
%% AJG-AWT:56
\begin{lemma}\label{le:positively-span-iff-I}
  Given a set $\{v_1,\dots,v_M\}$ of $M$ arbitrary vectors in $\real^N$,
  then the following statements are equivalent
  \begin{enumerate}
  \item $\{v_1,\dots,v_M\}$ positively spans $\real^N$;
  \item $ 0 \in \interior{\co \{v_1,\dots,v_M\}}$;
  \item for each $w \in \real^N$, there exists $v_i$ such that $w \cdot  v_i> 0$.
  \end{enumerate}
\end{lemma}

Let $Q$ be a convex polygon in $\real^2$. We denote by $\ED (Q) = \{
e_1,\dots,e_M \}$ and $\VE (Q) = \{ v_1,\dots,v_L\}$ the set of edges and
vertexes of $Q$, respectively.
%
% \sindex{Edges}{$\ED(Q)$}{Edges of polygon $Q$}%
% \sindex{Vertexes}{$\VE(Q)$}{Vertexes of polygon $Q$}%
%
Let $P = (p_1,\dots,p_n) \in Q^n\subset(\real^2)^n$ denote the location of
$n$ generators in the space $Q$. Let $\pi_i : Q^n \rightarrow Q$ be the
canonical projection onto the $i$th factor, $\pi_i (p_1,\dots,p_n) = p_i$.
%
%\sindex{Projection}{$\pi_i$}{Canonical projection from $Q^n$ onto the $i$th factor}%
%
Note that this mapping is surjective, continuous and open (the latter
meaning that open sets of $Q^n$ are mapped onto open sets of $Q$).

\subsection{Voronoi partitions}
We present here some relevant concepts on Voronoi diagrams and refer the
reader to~\cite{MdB-MvK-MO:97,AO-BB-KS-SNC:00} for comprehensive
treatments.  A \emph{partition} of $Q$ is a collection of $n$ polygons
$\WW=\{W_1,\dots,W_n\}$ with disjoint interiors whose union is $Q$.  Of
course, more general types of partitions could be considered (as, for
instance, continuous deformations of the previous ones), but these ones
will be sufficient for our purposes.  The \emph{Voronoi partition}
$\VV(P)=(V_1(P),\dots,V_n(P))$ of $Q$ generated by the points
$(p_1,\dots,p_n)$ is defined by:
\begin{equation*}
  V_i(P) = \setdef{q\in Q}{\|q-p_i\|\leq \|q-p_j\| \, , \; \forall j\neq i}.
\end{equation*}
%
%\sindex{Voronoi}{$\VV(P)$}{Voronoi partition of $Q$ generated by  $P=(p_1,\dots,p_n)$}
%
For simplicity, we shall refer to $V_i(P)$ as $V_i$.  Since $Q$ is a convex
polygon, the boundary of each $V_i$ is the union of a finite number of
segments.  If $V_i$ and $V_j$ share an edge, i.e., $V_i \cap V_j$ is
neither empty nor a singleton, then $p_i$ is called a \emph{(Voronoi)
  neighbor} of $p_j$ (and vice-versa).  All Voronoi neighboring relations
are encoded in the mapping $\NN : Q^n \times \until{n} \rightarrow
2^{\until{n}}$ where $\NN (P,i)$ is the set of indexes of the Voronoi
neighbors of $p_i$.  Of course, $j \in \NN (P,i)$ if and only if $i \in \NN
(P,j)$.  We will often omit $P$ and instead write $\NN (i)$.
%
%\sindex{NN}{$\NN (P,i)$, $\NN(i)$}{Set of neighbors of the $i$th
%  generator at configuration $P$}

\begin{figure}[htb] 
  \begin{center}
    \resizebox{.5\linewidth}{!}{\input{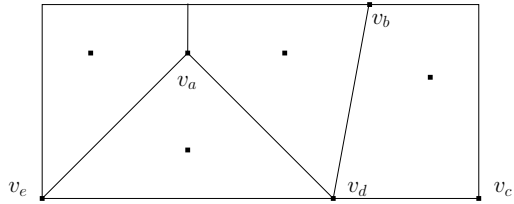}} 
  \end{center}
  \caption{A Voronoi partition with degenerate and nondegenerate
    vertexes. Vertexes $v_a$, $v_b$, and $v_c$ are nondegenerate vertexes
    of type (a), (b), (c), respectively. Vertexes $v_d$ and $v_e$ are
    degenerate.}
  \label{fig:voronoi-vertexes}
\end{figure}

For $P\in Q^n$, the vertexes of the Voronoi partition $\VV (P)$ are
classified as follows: the vertex $v$ is \emph{of type (a)} if it is the
center of the circle passing through three generators (say, $p_i$, $p_j$,
and $p_k$), the vertex $v$ is \emph{of type (b)} if it is the intersection
between an edge of $Q$ and the bisector determined by two generators (say,
$e$, $p_i$, and $p_j$), and the vertex $v$ is \emph{of type (c)} if it is a
vertex of $Q$, i.e., it is determined by two edges of $Q$ and by the
generator of a cell containing it (say, $e$, $f$, and $p_i$).
Correspondingly, we shall write $v(i,j,k)$, $v(e,i,j)$, and $v(e,f,i)$
respectively, whenever we are interested in making explicit the elements
defining the vertex $v$.
% %
% \sindex{v1}{$v(i,j,k)$}{Vertex of $\VV (P)$ determined by $p_i$, $p_j$ and
%   $p_k$}
% \sindex{v2}{$v(e,i,j)$}{Vertex of $\VV (P)$ determined by $e \in \ED (Q)$ and
%   $p_i$, $p_j$}
% \sindex{v3}{$v(e,f,i)$}{Vertex of $\VV (P)$ determined by $e,f \in \ED (Q)$
%   and $p_i$}
%
The vertex $v\in \VE (V_i (P))$ is said to be \emph{nondegenerate} if it is
determined by exactly three elements (e.g., as described above, either
three generators, or an edge and two generators, or two edges and one
generator), otherwise it is said to be \emph{degenerate}. Further, the
configuration $P$ is said to be \emph{nondegenerate at the $i$th generator}
if all vertexes $v\in\VE(V_i(P))$ are nondegenerate, otherwise $P$ is
\emph{degenerate at the $i$th generator}. Finally, a configuration $P$ is
said to be \emph{nondegenerate} if all its vertexes are nondegenerate,
otherwise it is said to be \emph{degenerate}. These concepts are
illustrated in Fig.~\ref{fig:voronoi-vertexes}.

For $P\in Q^n$, the edges of the Voronoi partition $\VV (P)$ are classified
as follows: the edge $e$ is \emph{of type (a)} if it is a segment of the
bisector determined by two generators (say, $p_i$, $p_j$), the edge $e$ is
\emph{of type (b)} if it is contained in the boundary of $Q$, i.e., it is a
subset of an edge of $Q$ and it belongs to a single cell (say, the cell of
the generator $p_i$).  Correspondingly, we shall write $e(i,j)$ and $e(i)$
respectively, whenever we are interested in making explicit the elements
defining the edge $e$.
%
% \sindex{eij}{$e(i,j)$}{Edge of $\VV (P)$ determined by $p_i$ and $p_j$}
% \sindex{ei}{$e(i)$}{Edge of $\VV (P)$ belonging to $V_i$ and to the boundary of
%   $Q$}
%
Further, when considering an edge of type (a), we let $n_{e(i,j)}$ denote
the unit normal to $e(i,j)$ pointing toward $\interior{V_i(P)}$. When
considering an edge of type (b), we let $n_{e(i)}$ denote the unit normal
to $e(i)$ pointing toward $\interior{Q}$.  
% %
% \sindex{Normal-i}{$n_{e(i,j)}$}{Unit normal to $e(i,j)$  pointing toward  $\interior{V_i(P)}$}
% \sindex{Normal-ii}{$n_{e(i)}$}{Unit normal to $e(i)$ pointing toward  $\interior{Q}$}
%

\subsection{The disk-covering and the sphere-packing problems}
\label{se:locational-functions}

We are interested in the following locational optimization problems
\begin{align}
  & \min_{p_1,\ldots,p_n} \left\{ \max_{q\in Q} \left\{ \min_{i\in\until{n}}
      \|q-p_i\| \right\}  \right\} \, , \label{eq:disk-covering} \\
  & \max_{p_1,\ldots,p_n} \left\{ \min_{
      \begin{subarray}{c}
        i,j\in\until{n}\\
        i \neq j , \, e\in\ED(Q)
      \end{subarray}
      } \left\{ {\textstyle \frac{1}{2}} \|p_i-p_j\| , \D_e(p_i)\right\}
  \right\} \, .
  \label{eq:sphere-packing} 
\end{align}
The optimization problem~\eqref{eq:disk-covering} is referred to as the
$p$-center problem in~\cite{ZD:95,AS-ZD:96}.  Throughout the paper, we will
refer to it as the multi-circumcenter problem.  In the context of coverage
control of mobile sensor networks~\cite{JC-SM-TK-FB:02j}, the
multi-circumcenter problem corresponds to considering the worst case
scenario, in which no information is available on the distribution of the
events taking place in the environment $Q$.  The network therefore tries to
minimize the largest possible distance of any point in $Q$ to one of the
generators' locations given by $p_1, \dots,p_n$, i.e.  to minimize the
function,
\[
\HHDC (P) = \max_{q\in Q} \left\{ \min_{i\in\until{n}} \|q-p_i\| \right\} =
\max_{i\in\until{n}} \left\{ \max_{q \in V_i} \|q-p_i\| \right\} \, .
\]
It is conjectured in~\cite{AS-ZD:96} that this problem can be restated as a
disk-covering problem: how to cover a region with (possibly overlapping)
disks of minimum radius. The disk-covering problem then reads:
\begin{equation*}
  \min \setdef{R}{\cup_{i \in \until{n}} \ov{B}_2(p_i,R) \supseteq Q} \, .
\end{equation*}
We shall present a proof of this statement in Theorem~\ref{th:minima-DC}
below. Given a polytope $W$ in $\real^N$, its circumcenter, denoted by $\CC
(W)$, is the center of the minimum-radius sphere that contains $W$. The
circumradius of $W$, denoted by $\CR (W)$, is the radius of this sphere.
%
%
% \sindex{HHDC}{$\HHDC$}{Multi-circumcenter function}
% \sindex{CC}{$\CC(Q)$}{Circumcenter of polytope $Q$}
% \sindex{CR}{$\CR(Q)$}{Circumradius of polytope $Q$}
% %
We will say that $P$ is a \emph{circumcenter Voronoi configuration} if $p_i =
\CC (V_i(P))$, for all $i \in \until{n}$. We denote by $\VE_{\DC}(\VV(P))$
the set of vertexes of the Voronoi partition where the value $\HHDC (P)$ is
attained, i.e. $v \in \VE_{\DC}(\VV(P))$ if there exists $i$ such that $v\in
V_i (P)$ and $\| v-p_i\| = \HHDC(P)$.
% %
% \sindex{VEDC}{$\VE_{\DC}(\VV(P))$}{Vertexes of $\VV (P)$ where the value of
%   $\HHDC (P)$ is attained}
% \sindex{HHSP}{$\HHSP$}{Multi-incenter function}
% %

We will refer to the optimization problem~\eqref{eq:sphere-packing} as the
multi-incenter problem. In the context of applications, this problem
corresponds to the situation where we are interested in maximizing the
coverage of the area $Q$ in such a way that the sensing radius of the
generators do not overlap (in order not to interfere with each other) or
leave the environment. We therefore consider the maximization of the function
\[
\HHSP (P) = \min_{
      \begin{subarray}{c}
        i,j\in\until{n}\\
        i \neq j , \, e\in\ED(Q)
      \end{subarray}
      } \left\{ {\textstyle \frac{1}{2}} \|p_i-p_j\| , \D_e(p_i)\right\} =
    \min_{i\in\until{n}} \left\{ \min_{ q \not \in \interior{V_i}} \|q-p_i\|
    \right\} \, .
\]
A similar conjecture to the one presented above is that the multi-incenter
problem can be restated as a sphere-packing problem: how to maximize the
coverage of a region with non-overlapping disks (contained in the region) of
minimum radius.  The problem reads:
\begin{equation*}
  \max \setdef{R}{\cup_{i \in \until{n}} \ov{B}_2(p_i,R) \subseteq  Q \, , \;
  B_2(p_i,R) \cap B_2(p_j,R) = \emptyset} \, . 
\end{equation*}
In Theorem~\ref{th:maxima-SP} we provide a positive answer to this question.
Given a polytope $W$ in $\real^N$, its incenter set (or Chebyshev center set,
see~\cite{SB-LV:02}), denoted by $\IC (W)$, is the set of the centers of
maximum-radius spheres contained in $W$. The inradius of $W$, denoted by $\IR
(W)$, is the common radius of these spheres.
%
% \sindex{IC}{$\IC(Q)$}{Incenter set of polytope $Q$}
% \sindex{IR}{$\IR(Q)$}{Inradius of polytope $Q$}
%
We will say that $P \in Q^n$ is an \emph{incenter Voronoi configuration} if
$p_i \in \IC (V_i(P))$, for all $i \in \until{n}$. If $P$ is an incenter
Voronoi configuration, and each Voronoi region $V_i(P)$ has a unique
incenter, $\IC (V_i(P)) =\{p_i\}$, then we will say that $P$ is a
\emph{generic incenter Voronoi configuration}. We denote by
$\ED_{\SP}(\VV(P))$ the set of edges of the Voronoi partition where the value
$\HHSP (P)$ is attained, i.e. $e \in \ED_{\SP}(\VV(P))$ if there exists $i$
such that $e \in \ED(V_i (P))$ and $\D_e (p_i) = \HHSP(P)$.
%
% \sindex{EDSP}{$\ED_{\SP}(\VV(P))$}{Edges where the value of $\HHSP(P)$ is
%   attained}
%

\subsection{Nonsmooth analysis}\label{sec:nonsmooth analysis}

The following facts on nonsmooth analysis~\cite{FHC:83} will be most helpful
in analyzing the properties of the locational optimization functions for the
disk-covering and the sphere-packing problems, as well as the convergence of
the distributed algorithms we will propose to extremize them.
\begin{definition}
  A function $f:\real^N \rightarrow \real$ is said to be \emph{locally
    Lipschitz near $x\in \real^N$} if there exist positive constants $L_x$
  and $\eps$ such that $| f(y) - f(y') | \le L_x \| y-y' \|$ for all $y,y'
  \in B_N(x,\eps)$.
\end{definition}

Note that continuously differentiable functions at $x$ are locally Lipschitz
near $x$. The usual \emph{right directional derivative} of $f$ at $x$ in the
direction of $v \in \real^N$ is defined as
\begin{align*}
  f'(x,v) = \lim_{t \rightarrow 0^+} \frac{f(x+tv) - f(x)}{t} \, ,
\end{align*}
when this limits exists. On the other hand, the \emph{generalized directional
  derivative} of $f$ at $x$ in the direction of $v \in \real^N$ is defined as
\begin{align*}
  f^o(x;v) = \limsup_{
    \begin{subarray}{l}
      y \rightarrow x\\
      t \rightarrow 0^+
    \end{subarray}
    } \frac{f(y+tv)-f(y)}{t} \, .
\end{align*}
This notion of directional derivative has the advantage of always being
well-defined.
\begin{definition}
  A function $f:\real^N \rightarrow \real$ is said to be \emph{regular at
    $x\in \real^N$} if for all $v \in \real^N$, $f'(x;v)$ exists and
  $f^o(x;v) = f'(x;v)$.
\end{definition}

Again, a continuously differentiable function at $x$ is regular at $x$. Also,
a locally Lipschitz function at $x$ which is convex (or concave) is also
regular (cf. Proposition~2.3.6 in~\cite{FHC:83}).

From Rademacher's Theorem~\cite{FHC:83}, we know that locally Lipschitz
functions are continuously differentiable almost everywhere (in the sense of
Lebesgue measure). If $\Omega_f$ denotes the set of points in $\real^N$ at
which $f$ fails to be differentiable, and $S$ denotes any other set of
measure zero, the \emph{generalized gradient} of $f$ is defined by
\[
\partial f (x) = \co \setdef{\lim_{i \rightarrow +\infty} df (x_i)}{x_i
  \rightarrow x \, , \; x_i \not \in S \cup \Omega_f} \, .
\]
%
% \sindex{Generalizedgradient}{$\partial f$}{Generalized gradient of the
%   locally Lipschitz function $f$}%
%
Note that this definition coincides with $df (x)$ if $f$ is continuously
differentiable at $x$. The generalized gradient and the generalized
directional derivative (cf.  Proposition~2.1.2 in~\cite{FHC:83}) are related
by $f^o (x;v) = \max \setdef{\zeta \cdot v}{\zeta \in \partial f (x)}$, for
each $v \in \real^N$. A point $x \in\real^N$ which verifies that $ 0 \in
\partial f(x)$ is called a \emph{critical point of $f$}.

The following result corresponds to Proposition~2.3.12 in~\cite{FHC:83}.
\begin{proposition}\label{prop:gradient-F-nice}
  Let $\setdef{f_k:\real^N \rightarrow \real}{k \in \until{m}}$ be a finite
  collection of locally Lipschitz functions near $x \in \real^N$.  Consider
  $f(x') = \min \setdef{f_k(x')}{k \in \until{m}}$.  Then,
  \begin{enumerate}
  \item $f$ is locally Lipschitz near $x$,
  \item if $I(x')$ denotes the set of indexes $k$ for which $f_k(x')=f(x')$,
    we have,
    \begin{align}\label{eq:gradient-F-nice}
      \partial f (x) \subset \co \setdef{\partial f_i (x)}{i \in I(x)} \, ,
    \end{align}
    and if each $f_i$ is regular at $x$ for $i \in I(x)$, then equality holds
    and $f$ is regular at $x$.
  \end{enumerate}
\end{proposition}

The extrema of Lipschitz functions are characterized by the following result.
\begin{proposition}\label{prop:zero-in-gradient}
  Let $f$ be a locally Lipschitz function at $x \in \real^N$.  If $f$ attains
  a local minimum or maximum at $x$, then $0 \in \partial f (x)$, i.e., $x$
  is a critical point.
\end{proposition}

Let $\LN: 2^{\real^N} \rightarrow 2^{\real^N}$ be the set-valued mapping
that associates to each subset $S$ of $\real^N$ the set of its least-norm
elements $\LN (S)$.  If the set $S$ is convex, then the set $\LN (S)$
reduces to a singleton and we note the equivalence $\LN(S)=\proj{S}(0)$.
Along the paper, we shall only apply this function to convex sets. For a
locally Lipschitz function $f$, we consider the \emph{generalized gradient
  vector field} $\LN (\partial f) : \real^N \rightarrow \real^N$ given by
$x\mapsto\LN (\partial f) (x) = \LN (\partial f (x))$.
%
%\sindex{Leastnorm}{{$\LN (S)$}}{Least-norm element of the convex set $S$}%
%

\begin{theorem}\label{th:descent-direction}
  Let $f$ be a locally Lipschitz function at $x$. Assume $0 \not \in \partial
  f (x)$.  Then, there exists $T>0$ such that
  \begin{equation*}
    f(x-t \, \LN (\partial f)(x)) \le f(x) - \frac{t}{2} \|\LN (\partial f)(x)\|^2 \, , \quad 0
    < t <  T \, . 
  \end{equation*}
  The vector $-\LN (\partial f)(x)$ is called a \emph{direction of descent}.
\end{theorem}

\subsection{Stability analysis via nonsmooth Lyapunov functions}\label{se:stability-via-Lipschitz} 

For differential equations with discontinuous right-hand sides, solutions
are defined in terms of differential inclusions~\cite{AFF:88}.

Let $F: \real^N \rightarrow 2^{\real^N}$ be a set-valued map. Consider the
differential inclusion
\begin{align}\label{eq:diff-inclusion}
  \dot{x} \in F (x) \, .
\end{align}
A solution to this equation on an interval $[t_0,t_1] \subset \real$ is
defined as an absolutely continuous function $x:[t_0,t_1] \rightarrow
\real^N$ such that $\dot{x}(t) \in F(x(t))$ for almost all $t \in [t_0,t_1]$.
Given $x_0 \in \real^N$, the existence of at least a solution with initial
condition $x_0$ is guaranteed by the following lemma.

\begin{lemma}\label{le:existence-solution}
  Let the mapping $F$ be upper semicontinuous with nonempty, compact and
  convex values. Then, given $x_0 \in \real^N$, there exists at least a
  solution of~\eqref{eq:diff-inclusion} with initial condition $x_0$.
\end{lemma}

Now, consider the differential equation
\begin{align}\label{eq:diff-equation}
  \dot{x}(t) = X(x(t)) \, ,
\end{align}
where $X:\real^N \rightarrow \real^N$ is measurable and essentially locally
bounded. The solution of this equation has to be understood in the Filippov
sense. For each $x \in \real^N$, consider the set
\begin{align*}
  K[X](x) = \bigcap_{\delta >0} \bigcap_{\mu (S) = 0} \co \{
  X(B_N(x,\delta)\setminus S) \} \, ,
\end{align*}
%
% \sindex{K}{{$K[X]$}}{Filippov mapping associated with a measurable and
%   essentially locally bounded mapping $X:\real^N \rightarrow \real^N$ }%
%
where $\mu$ denotes the usual Lebesgue measure in $\real^N$. Alternatively,
one can show~\cite{BP-SSS:87} that there exists a set $S_X$ of measure zero
such that
\begin{align*}
  K[X](x) = \co \setdef{\lim_{i \rightarrow +\infty} X(x_i)}{x_i
    \rightarrow x \, , \; x_i \not \in S \cup S_X} \, ,
\end{align*}
where $S$ is any set of measure zero. A Filippov solution
of~\eqref{eq:diff-equation} on an interval $[t_0,t_1] \subset \real$ is
defined as a solution of the differential inclusion
\begin{align}\label{eq:Filippov-diff-inclusion}
  \dot{x} \in K[X](x) \, .
\end{align}
Since the multivalued mapping $K[X]: \real^N \rightarrow 2^{\real^N}$ is
upper semicontinuous with nonempty, compact, convex values and locally
bounded (cf.~\cite{AFF:88}), the existence of Filippov solutions
of~\eqref{eq:diff-equation} is guaranteed by
Lemma~\ref{le:existence-solution}.

A set $M$ is \emph{weakly invariant} (respectively \emph{strongly invariant})
for~\eqref{eq:diff-equation} if for each $x_0 \in M$, $M$ contains a maximal
solution (respectively all maximal solutions) of~\eqref{eq:diff-equation}.
Given a locally Lipschitz function $f:\real^N \rightarrow \real$, the
\emph{set-valued Lie derivative of $f$ with respect to $X$} at $x$ is defined
as
\begin{align*}
  \setLieder{X}{f} (x) = \setdef{a \in \real}{\exists v \in K[X](x) \;
    \text{such that} \; \zeta \cdot v = a \, , \; \forall \zeta \in \partial
    f(x)} \, .
\end{align*}
%
% \sindex{setLieder}{$\setLieder{X}{f}$}{Set-valued Lie derivative of $f$ with respect to $X$}%
%
For each $x \in \real^N$, $\setLieder{X}{f}(x)$ is a closed and bounded
interval in $\real$, possibly empty. If $f$ is continuously differentiable at
$x$, then $\setLieder{X}{f}(x) = \setdef{df \cdot v}{v \in K[X](x)}$. If, in
addition, $X$ is continuous at $x$, then $\setLieder{X}{f}(x)$ corresponds to
the singleton $\{ \Lc_{X} f (x) \}$, the usual Lie derivative of $f$ in the
direction of $X$ at $x$. The importance of the set-valued Lie derivative
stems from the next result~\cite{AB-FC:99}.

\begin{theorem}\label{th:single-trajectory}
  Let $x:[t_0,t_1] \rightarrow \real^N$ be a Filippov solution
  of~\eqref{eq:diff-equation}. Let $f$ be a locally Lipschitz and regular
  function. Then $\frac{d}{dt} \left( f (x(t))\right)$ exists a.e. and
  $\frac{d}{dt} \left( f (x(t))\right) \in \setLieder{X}{f} (x(t))$~a.e.
\end{theorem}

The following result is a generalization of LaSalle
principle for differential equations of the form~\eqref{eq:diff-equation} with
nonsmooth Lyapunov functions. The formulation is taken from~\cite{AB-FC:99},
and slightly generalizes the one presented in~\cite{DS-BP:94}.

\begin{theorem}[LaSalle principle] \label{th:LaSalle}
  Let $\map{f}{\real^N}{\real}$ be a locally Lipschitz and regular
  function.  Let $x_0 \in \real^N$ and let $f^{-1} (\le f(x_0),x_0)$ be the
  connected component of $\setdef{x \in \real^N}{f(x) \le f(x_0)}$
  containing $x_0$. Assume the set $f^{-1} (\le f(x_0),x_0)$ is bounded and
  assume either $\max \setLieder{X}{f}(x) \le 0$ or $\setLieder{X}{f}(x) =
  \emptyset$ for all $x \in f^{-1} (\le f(x_0),x_0)$. Then $f^{-1} (\le
  f(x_0),x_0)$ is strongly invariant for~\eqref{eq:diff-equation}.  Let
  \[
  Z_{X,f} = \setdef{x \in \real^N}{0 \in \setLieder{X}{f}(x)} \, .
  \]
%   %                              
%   \sindex{Z}{$Z_{X,f}$}{Set formed by points $x \in\real^N$ such that  $0$
%     belongs to  $\setLieder{X}{f}(x)$}%
%   %
  Then, any solution $x:[t_0,+\infty) \rightarrow \real^N$
  of~\eqref{eq:diff-equation} starting from $x_0$ converges to the largest
  weakly invariant set $M$ contained in $\ov{Z}_{X,f} \cap f^{-1} (\le
  f(x_0),x_0)$. Furthermore, if the set $M$ is a finite collection of points,
  then the limit of all solutions starting at $x_0$ exists and equals one of
  them.
\end{theorem}

The proof of the last fact in the theorem statement is the same as in the
smooth case, since it only relies on the continuity of the trajectory. The
next statement is based on Theorem~2 of~\cite{BP-SSS:87}.

\begin{proposition}\label{prop:finite-time} 
  Under the same assumptions of Theorem~\ref{th:LaSalle}, if $\max
  \setLieder{X}{f}(x) < -\eps <0 $ a.e. on $\real^N \setminus Z_{X,f}$,
  then $Z_{X,f}$ is attained in finite time.
\end{proposition}
\begin{proof}
  Let $x:[t_0,+\infty) \rightarrow \real^N$ be a Filippov solution starting
  from $x_0$. We argue that there must exist $T$ such that $x(T) \in
  Z_{X,f}$.  Otherwise, we have
  \[
  f (x (t)) = f(x(t_0)) + \int_{t_0}^t \frac{d}{ds} f (x (s)) ds < f
  (x(t_0))- \eps (t- t_0) \stackrel{t \rightarrow +\infty}{\longrightarrow}
  -\infty \, ,
  \]
  which contradicts the fact that $f^{-1} (\le f(x_0),x_0)$ is strongly
  invariant and bounded.  \quad
\end{proof}

\subsection{Nonsmooth gradient flows}

Finally, we are in a position to present the nonsmooth analogue of well-known
results on gradient flows.  Given a locally Lipschitz and regular function
$f$, consider the following generalized gradient flow
\begin{equation}\label{eq:natural-gradient}
  \dot{x}(t) = - \LN (\partial f)(x(t)) \, . 
\end{equation}
Theorem~\ref{th:descent-direction} guarantees that, unless the flow is at a
critical point, $-\LN (\partial f)(x)$ is always a direction of descent at
$x$. In general, the vector field $\LN (\partial f)$
in~\eqref{eq:natural-gradient} is discontinuous, and therefore its solution
must be understood in the Filippov sense. Note that, since $f$ is locally
Lipschitz, $\LN (\partial f) = df$ almost everywhere. An important
observation in this setting is that $K[df](x) = \partial f (x)$
(cf.~\cite{BP-SSS:87}). The following result, which is a generalization of
the discussion in~\cite{AB-FC:99}, guarantees the convergence of this flow to
the set of critical points of $f$.

\begin{proposition}\label{prop:natural-gradient}
  Let $x_0 \in \real^N$ and assume $f^{-1} (\le f(x_0),x_0)$ is bounded.
  Then, any solution $x:[t_0,+\infty) \rightarrow \real^N$ of
  eq.~\eqref{eq:natural-gradient} starting from $x_0$ converges
  asymptotically to the set of critical points of $f$ contained in $f^{-1}
  (\le f(x_0),x_0)$.
\end{proposition}

\begin{proof}
  Let $a \in \setLieder{-\LN (\partial f)}{f}(x)$. By definition, there
  exists $w \in K[-\LN (\partial f)](x) = -\partial f (x)$ such that $a=w
  \cdot \zeta$ for all $\zeta \in \partial f (x)$. In particular, for
  $\zeta=-w \in \partial f (x)$, we have $a = - \| w \|^2 \le 0$. Therefore,
  $\max \setLieder{-\LN (\partial f)}{f}(x) \le 0$ or $\setLieder{-\LN (\partial
    f)}{f}(x) = \emptyset$. Now, resorting to the LaSalle principle
  (Theorem~\ref{th:LaSalle}), we deduce that any solution $x:[t_0,+\infty)
  \rightarrow \real^N$ starting from $x_0$ converges to the largest weakly
  invariant set contained in $\ov{Z}_{-\LN (\partial f),f} \cap f^{-1} (\le
  f(x_0),x_0)$. Let us see that $Z_{-\LN (\partial f),f}$ is equal to $L_0 =
  \setdef{x \in Q^n}{0 \in \partial f(x)}$.  Obviously, $L_0 \subset Z_{-\LN
    (\partial f),f}$.  Conversely, assume $x \in Z_{-\LN (\partial f),f}$.
  Then, $0 \in \setLieder{-\LN (\partial f)}{f}(x)$, i.e., there exists $v
  \in - \partial f (x)$ such that $\zeta \cdot v = 0$ for all $\zeta \in
  \partial f (x)$. In particular, for $\zeta = - v$, we get $\|v\|^2 = 0$,
  that is, $v = 0 \in \partial f (x)$, as desired.  Note that $Z_{-\LN
    (\partial f),f}=L_0$ is the equilibrium set
  of~\eqref{eq:natural-gradient} and therefore is weakly invariant.  Finally,
  we prove that it is also closed. Let $x \in \ov{Z}_{-\LN (\partial f),f}$
  and consider a sequence $\setdef{x_k \in \real^N}{k \in \natural} \subset
  Z_{-\LN (\partial f),f}$ such that $x_k \rightarrow x$.  Then, using the
  fact that the multivalued mapping $K[-v]$ is upper semicontinuous, for any
  $\eps>0$, there exists $k_0$ such that for $k \ge k_0$, $\partial f (x_k)
  \subset \partial f (x) + B_N(0,\eps)$.  Since $x_k \in Z_{-\LN (\partial
    f),f}$, then $0 \in \partial f (x) + B_N(0,\eps)$ for all $\eps >0$, and
  this implies that $0 \in \partial f (x)$, i.e., $x \in Z_{-\LN (\partial
    f),f}$.  Hence the largest weakly invariant set contained in
  $\ov{Z}_{-\LN (\partial f),f} \cap f^{-1} (\le f(x_0),x_0)$ is $Z_{-\LN
    (\partial f),f} \cap f^{-1} (\le f(x_0),x_0) = \setdef{x \in f^{-1} (\le
    f(x_0),x_0)}{0 \in \partial f (x)}$.  \quad
\end{proof}

\section{The 1-center problems}\label{sec:1-center}

In this section we consider the disk-covering and the sphere-packing
problems with a single generator, i.e., $n=1$. This treatment will give us
the necessary insight to tackle later the more involved multi-center
version of both problems. When $n=1$, the minimization of $\HHDC$ simply
consists of finding the center of the minimum-radius sphere enclosing the
polygon $Q$. On the other hand, the maximization of $\HHSP$ corresponds to
determining the center of the maximum-radius sphere contained in $Q$. Let
us therefore define the functions
\begin{align}
  \lg_Q (p) &= \max \setdef{\| q - p \|}{q \in Q} = \max \setdef{\| v - p
    \|}{v \in \VE (Q)} \, , \nonumber \\
  \sm_Q (p) &= \min \setdef{\| q - p \|}{q \not \in \interior{Q}} = \min
  \setdef{ \D_e (p) }{e \in \ED (Q) } \, . \label{eq:1-center-alternative}
\end{align}
When $n=1$, we then have that $\HHDC = \lg_Q: Q \rightarrow \real$ and $\HHSP
= \sm_Q : Q \rightarrow \real$.
%
% \sindex{lg}{$\lg_Q(p)$}{Largest distance from $p$ to the boundary of
%   $Q$}
% \sindex{sm}{$\sm_Q(p)$}{Smallest distance from $p$ to the boundary of
%   $Q$}
% %

\subsection{Smoothness and critical points}
We here discuss the smoothness properties and the critical points of the
1-center functions.  Since the function $\lg_Q$ is the maximum of a
(finite) set of convex functions in $p$, it is also a convex
function~\cite{SB-LV:02}. Therefore, any local minimum of $\lg_Q$ is also
global.

\begin{lemma}
  The function $\lg_Q$ has a unique global minimum, which is the
  circumcenter of the polygon $Q$.
\end{lemma}

\begin{proof}
  Let $F : \real \rightarrow \real$ be any continuous non-decreasing
  function. Then,
  \[
  F (\lg_Q(p)) = \max \setdef{F(\| v - p \|)}{v \in \VE (Q) } \, .
  \]
  If we take $F(x)=x^2$, each function $\| v - p \|^2$ is strictly convex,
  and hence $F (\lg_Q(p))$ is also strictly convex.  Therefore, this latter
  function has a single minimum on $Q$. Since any global minimum of $\lg_Q$
  is also a global minimum of $F (\lg_Q(p))$, we conclude the result. \quad
\end{proof}

The function $\sm_Q$ is the minimum of a (finite) set of affine (hence,
concave) functions defined on the half-planes determined by the edges of $Q$,
and hence it is also a concave function~\cite{SB-LV:02} on the intersection
of their domains, which is precisely $Q$. Therefore, any local maximum of
$\sm_Q$ is also global. However, this maximum is not unique in general.

\begin{lemma}\label{le:incenter-set-convex-segment}
  The incenter set of the polygon $Q$ is the set of maxima of the function
  $\sm_Q$ and it is a segment.
\end{lemma}

\begin{proof}
  It is clear that the set of maxima of $\sm_Q$ is $\IC (Q)$. As a
  consequence of the concavity of $\sm_Q$ over the convex domain $Q$, one
  deduces that $\IC (Q)$ is a convex set.
  % Let $p_1$, $p_2 \in \IC (Q)$. For each $\lambda \in [0,1]$, we have
%   $\sm_Q (\lambda p_1 + (1-\lambda) p_2) \le \sm_Q (p_1) =
%   \sm_Q (p_2)$.  Since $\sm_Q$ is concave, we also have
%   $\sm_Q (\lambda p_1 + (1-\lambda) p_2) \ge \lambda \sm_Q
%   (p_1) + (1-\lambda) \sm_Q (p_2) = \sm_Q (p_1)$.
%   Therefore, $\sm_Q (\lambda p_1 + (1-\lambda) p_2) =
%   \sm_Q (p_1)$, and $\lambda p_1 + (1-\lambda) p_2 \in \IC(Q)$.
%   Hence, $\IC (Q)$ is a convex set.
  Now, assume there are three points $p_1, p_2, p_3$ in $\IC (Q)$ which are
  not aligned. Since $B_2 (q,\IR (Q)) \subset Q$ for all $q \in \co
  (p_1,p_2,p_3) \subset \IC (Q)$, and $\co (p_1,p_2,p_3)$ has non-empty
  interior, there exist $q_0 \in Q$ and $r > \IR(Q)$ such that $B_2(q_0,r)
  \subset Q$, which is a contradiction. \quad
\end{proof}

Note that the circumcenter of a polygon can be computed via the finite-step
algorithm described in~\cite{SS:91}.  The incenter set of a polygon can be
computed via the following linear program: maximize the inradius subject to
the constraints that the distance between the incenter and each of the
polygon edges must be greater than or equal to the inradius. In what
follows, let us examine dynamical systems that compute these geometric
centers.

\begin{proposition}\label{prop:remarkable}
  The functions $\lg_Q(p)$, $\sm_Q(p)$ are locally Lipschitz and regular, and
  their generalized gradients are given by
  \begin{align}
    \partial \lg_Q (p) &= \co \setdef{\versus{p-v}}{v \in \VE (Q) \, , \;
      \lg_Q (p)= \| p - v \|} \, ,
    \label{eq:partial-HHDC-1-center} \\
    \partial \sm_Q (p) &= \co \setdef{n_e}{e \in \ED (Q) \, , \; \sm_Q (p)=
      \D_e (p)} \, .
    \label{eq:partial-HHSP-1-center}
  \end{align}
  Moreover, 
  \begin{gather}\label{eq:remarkable}
    0 \in \partial \lg_Q(p) \Longleftrightarrow p = \CC (Q) \, , \quad 0 \in
    \partial \sm_Q(p) \Longleftrightarrow p \in \IC (Q) \, ,
  \end{gather}
  and, if $0 \in \interior{\partial \sm_Q(p)}$, then $\IC (Q) = \{p\}$.
\end{proposition}

\begin{proof}
  Given the expressions in~\eqref{eq:1-center-alternative} and
  Proposition~\ref{prop:gradient-F-nice}, we deduce that $\lg_Q$ and $\sm_Q$
  are locally Lipschitz and regular, and that their generalized gradients are
  respectively given by~\eqref{eq:partial-HHDC-1-center}
  and~\eqref{eq:partial-HHSP-1-center}.  Concerning~\eqref{eq:remarkable},
  the implications from right to left in~\eqref{eq:remarkable} readily follow
  from Proposition~\ref{prop:zero-in-gradient}. As for the other ones, note
  that it is sufficient to prove that $p$ is a local minimum, respectively
  that $p$ is a local maximum. We prove the result for the function $\lg_Q$.
  The proof for $\sm_Q$ is analogous.  Assume that $0 \in \partial \lg_Q(p)$.
  Then, there exist vertexes $v_{i_1},\dots,v_{i_K}$ of $Q$ with $\lg_Q(p) =
  \| v_{i_l} - p \|$, $l \in \until{K}$ such that $0 = \sum_{l \in \until{K}}
  \lambda_l \versus{p-v_{i_l}}$, where $\sum_{l \in \until{K}} \lambda_l =
  1$, $\lambda_l \ge 0$, $l \in \until{K}$. Let $U$ be a neighborhood of $p$
  and take $q \in U$. One can show that there must exist $l^*$ such that
  $(p-v_{i_{l^*}}) \cdot (q-p) \ge 0$, since otherwise $0 = 0 \cdot (q-p) =
  (\sum_{l \in \until{K}} \lambda_l \versus{p-v_{i_l}}) \cdot (q-p) <0$,
  which is a contradiction. Then,
  \[
  \| q - v_{i_{l^*}} \|^2 = \| q - p \|^2 + \| p - v_{i_{l^*}} \|^2 -2
  (q-p)\cdot(v_{i_{l^*}}-p) \ge \| p - v_{i_{l^*}} \|^2 \, .
  \]
  Therefore, $ \lg_Q(q) \ge \| p - v_{i_{l^*}} \| = \lg_Q(p)$, which shows
  that $p$ is a local minimum.  Finally, if $0 \in \interior{\partial
    \sm_Q(p)}$, then one can see that $p$ is a strict local maximum.
  Furthermore, there cannot be any other local (hence global) maximum of
  $\sm_Q$, as we now show: assume $\bar{p} \in \IC (Q)$. By hypothesis, the
  sphere $B_2(\bar{p},\sm_Q(p))$ centered at $\bar{p}$ of radius $\sm_Q(p)$
  is contained in $Q$. Consider the vector $\bar{p} - p$. By
  Lemma~\ref{le:positively-span-iff-I}, there exists $e \in \ED (Q)$ with
  $\D_e (p) = \sm_Q(p)$ such that $(\bar{p} - p) \cdot n_e > 0$.
  Therefore, there are points of $B_2(\bar{p},\sm_Q(p))$ which necessarily
  belong to the half-plane defined by $e$ where $Q$ is not contained, which
  is a contradiction.  \quad
\end{proof}

\subsection{Convergence properties for nonsmooth gradient flows}

Here we study the generalized gradient flows arising from the two 1-center
functions.  An immediate consequence of
Propositions~\ref{prop:natural-gradient} and~\ref{prop:remarkable} is the
following result.

\begin{corollary}\label{coro:convergence-1-center}
  The gradient flows of the functions $\lg_Q$ and $\sm_Q$
  \begin{align}
    \dot{x}(t) & = - \LN (\partial \lg_Q)(x(t)) \, , \label{eq:natural-gradient-DC} \\
    \dot{x}(t) & = \LN (\partial \sm_Q)(x(t)) \, ,
    \label{eq:natural-gradient-SP}
  \end{align}
  converge asymptotically to the circumcenter $\CC (Q)$ and the incenter set
  $\IC (Q)$, respectively.
\end{corollary}

The following two propositions discuss the convergence properties of the
gradient descents.

\begin{proposition}\label{prop:finite-time-DC}
  If $0 \in \interior{\partial \lg_Q (\CC(Q))}$, then the
  flow~\eqref{eq:natural-gradient-DC} reaches $\CC (Q)$ in finite time.
\end{proposition}

\begin{proof}
  Let us prove that there exists $\eps > 0$ such that $\max
  \setLieder{-\LN[\lg_Q]}{\lg_Q} < -\eps$\/ a.e. on $Q \setminus
  \{ \CC (Q) \}$. Take $p \not = \CC (Q)$. We know that each element $a \in
  \setLieder{-\LN[\lg_Q]}{\lg_Q} (p)$ can be expressed as $a = -
  \| w \|^2$, with $-w \in \partial \lg_Q (p)$. Therefore, we have
  \[
  \max \setLieder{-\LN[\lg_Q]}{\lg_Q} (p) = - \| \LN[\lg_Q]
  (p) \|^2 \, .
  \]
  If there is a single vertex of $Q$ involved in $\partial \lg_Q (p)$,
  then moving along the direction $-\LN[\lg_Q] (p)$ obviously decreases
  the distance to that vertex while maintaining constant the norm of the
  least-norm element, which is $1$. If there are two or more vertexes
  involved, from the expression for the generalized gradient at $p$ (cf.
  eq.~\eqref{eq:partial-HHDC-1-center}), it is clear that one can express it
  as
  \[
  \setdef{x \in \real^N}{g_1(x) \le 0, \dots,g_s(x) \le 0} \, ,
  \]
  for some linear functions $g_r$. Note that the points $x \in \partial
  \lg_Q (p)$ such that $g_r (x) =0$ for some $r$ correspond to a set of
  the form $\co \{ \versus{p - v_{r,1}}, \versus{p-v_{r,2}} \}$, for certain
  vertexes $v_{r,1}$, $v_{r,2}$ of $Q$. Now, the computation of the least-norm
  element in $\partial \lg_Q (p)$ can be formulated as the convex
  problem,
  \begin{align*}
    & \text{minimize $\| x \|^2$}\\
    & \text{subject to $g_1(x) \le 0, \dots,g_s(x) \le 0$} \, .
  \end{align*}
  Let $x^* = \LN[\lg_Q] (p)$. Let $R$ denote the set of indexes $r$ for which
  $g_r(x^*)=0$. Then $x^*$ is a regular point~\cite{DGL:84}, meaning to say
  that $dg_r (x^*)$, $r \in R$ are linearly independent vectors. This is
  because the cardinality of $R$ is at most $2$ (since the intersection of
  two lines already determines a point), and the gradients of the functions
  $g_r$ are independent when considered pairwise.  We apply then the
  Kuhn-Tucker first-order necessary conditions for optimality~\cite{DGL:84}
  to conclude that there must exist $r^* \in R$ such that $g_{r^*}(x^*) = 0$.
  It is easy to see that $r^*$ must be unique, since otherwise $x^*$ does not
  have minimum norm. Therefore, we have that $\LN[\lg_Q] (p)$ is determined
  as the least-norm element in $\co \{ \versus{p - v_{r^*,1}},
  \versus{p-v_{r^*,2}} \}$. As a consequence, moving along the direction $-
  \LN[\lg_Q] (p)$ decreases the distance to the vertexes $v_{r^*,1}$,
  $v_{r^*,2}$, and hence the norm of the least-norm element decreases.  If,
  along the flow~\eqref{eq:natural-gradient-DC}, a new vertex of $Q$ enters in
  the computation of $\partial \lg_Q(p(t))$, then there can be a jump in the
  norm of $\LN[\lg_Q] (p(t))$, which by definition will always be decreasing.
  Finally, note that if $v_{r^*,1}$, $v_{r^*,2}$ are active at the
  circumcenter, then they cannot be opposite with respect to $\CC (Q)$
  precisely because of the assumption that $0$ lies in $\interior{\partial
  \lg_Q (\CC(Q))}$.  Therefore, we conclude
  \begin{multline*}
    \| \LN[\lg_Q] (p) \| \ge \eps =\min \big \{ 1, \{ \| \LN (\co \{
    \versus{\CC(Q)-v}, \versus{\CC(Q)-w} \}) \| \; | \\
    v, w \in I(\CC(Q)), \CC(Q)-v \not = - (\CC(Q)-w) \} \big\} > 0 \, , \quad
    \forall p \not = \CC (Q) \, .
  \end{multline*}
  Resorting now to Proposition~\ref{prop:finite-time}, we deduce that the
  circumcenter $\CC (Q)$ is attained in finite time.  \quad
\end{proof}

\begin{remark}
  {\rm Note that if $0 \in \partial \lg_Q (\CC(Q)) \setminus \interior{\partial
      \lg_Q (\CC (Q))}$, then generically convergence is achieved over an
      infinite time horizon.}
\end{remark}

\begin{proposition}\label{prop:finite-time-SP}
  The flow~\eqref{eq:natural-gradient-SP} reaches the set $\IC(Q)$ in finite
  time.
\end{proposition}

\begin{proof}
  Let $p \not \in \IC (Q)$. We know $ \min \setLieder{\LN[\sm_Q]}{\sm_Q} (p)
  = \| \LN[\sm_Q] (p) \|^2$. Moreover, for all $p \not \in \IC (Q)$, we have
  \begin{gather*}
    \| \LN[\sm_Q] (p) \| \ge \eps =\min \left \{ 1, \setdef{ \| \LN (\co \{
        n_e , n_f \} )\|}{e, f \in \ED (Q), n_e \not = - n_f} \! \right\} \!
    > 0 .
  \end{gather*}
  Resorting to Proposition~\ref{prop:finite-time}, we deduce the desired
  result.  \quad
\end{proof}

\begin{figure}[htb]
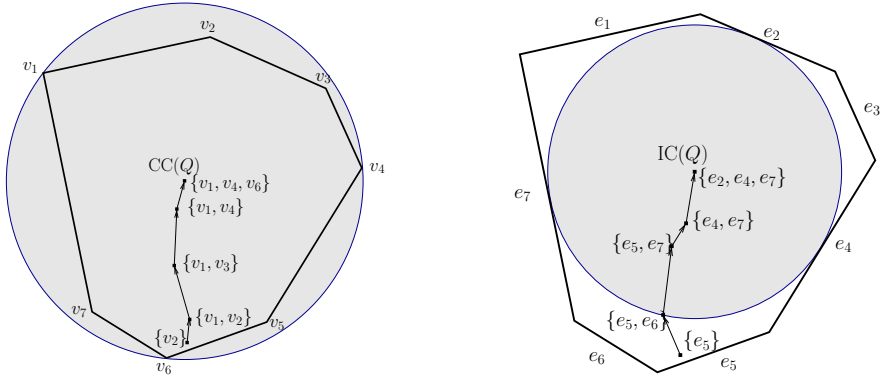
 
  \begin{center}
    \resizebox{.375\linewidth}{!}{\input{1-center-DC.tex}} \qquad\qquad\quad
    \resizebox{.375\linewidth}{!}{\input{1-center-SP.tex}}
  \end{center}
  \caption{Illustration of the gradient descent of $\lg_Q$ and $\sm_Q$. The
    points where the curve $t \mapsto p(t)$ fails to be differentiable
    correspond to points where there is a new vertex $v$ of $Q$ such that
    $\|p(t) -v\| = \lg_Q(p(t))$ (respectively a new edge $e$ of $Q$ such
    that $\D_e (p(t)) = \sm_Q(p(t))$). The circumcenter and the incenter
    are attained in finite time according to
    Propositions~\ref{prop:finite-time-DC}
    and~\ref{prop:finite-time-SP}.}\label{fig:center-gradient-descent}
\end{figure}

Fig.~\ref{fig:center-gradient-descent} shows an example of the
implementation of the gradient descent~\eqref{eq:natural-gradient-DC}
and~\eqref{eq:natural-gradient-SP}. Note that if the circumcenter $\CC (Q)$
(respectively the incenter set $\IC (Q)$) is first computed offline, then the
strategy of directly going toward it would converge in a less ``erratic" way.
Note also that the move-toward-the-center strategy is exponentially fast.

Finally, we conclude this section with four results useful for later
developments.
\begin{lemma}\label{le:simple-facts-DC-SP}
  Let $q \in Q$, let $v(q)$ be one of the vertexes of $Q$ which is furthest
  away from $q$, and let $e(q)$ be one of the edges of $Q$ which is nearest
  to $q$. Then,
  \begin{enumerate}
  \item $\LN[\lg_Q](q) \cdot (q-v(q)) \ge 0$, and the inequality is strict
    if $q \not = \CC (Q)$, 
  \item $(q-\CC(Q)) \cdot (q-v(q)) \ge \| q - \CC (Q) \|^2/2$,
  \item $\LN[\sm_Q](q) \cdot n_e \ge 0$, and the inequality is strict if $q
    \not \in \IC (Q)$, and
  \item $(x-q) \cdot n_e \ge \IR(Q) - \D_e(q) \ge 0$ for any $x \in \IC(Q)$,
    and the second inequality is strict if $q \not \in \IC(Q)$.
  \end{enumerate}
\end{lemma}
\begin{proof}
  Let $q$ be a point in $Q$.  If $q = \CC (Q)$, claims (i) and (ii) are
  obviously satisfied since $\LN[\lg_Q](q)=0$.  Assume then that $q \not =
  \CC (Q)$.  Let us prove first (i). By Proposition~\ref{prop:remarkable},
  $0 \not \in \partial \lg_Q (q)$, and hence $\LN[\lg_Q](q) \neq 0$.  Let
  us prove $\LN[\lg_Q](q) \cdot (q-v(q)) > 0$ reasoning by contradiction.
  If $\LN[\lg_Q](q) \cdot (q-v(q)) \le 0$, then $d/dt \left( \| q - t
    \LN[\lg_Q](q) - v \| \right)_{t=0} = \versus{q-v} \cdot
  (-\LN[\lg_Q](q)) \ge 0$, which implies that $\| q - t \LN[\lg_Q](q) - v
  \| \ge \| q - v \| = \lg_Q (q)$ for $t>0$ small enough. On the other
  hand, invoking Theorem~\ref{th:descent-direction}, we have that $\lg_Q
  (q) - t \| \LN[\lg_Q](q) \|^2/2 \ge \lg_Q (q - t \LN[\lg_Q](q)) \ge \|
  q-t\LN[\lg_Q](q) -v\|$.  Gathering both facts, we conclude $-t \|
  \LN[\lg_Q](q)\|^2/2 \ge 0$, which is a contradiction.
  
  Let us now prove (ii). Since $q \not = \CC (Q)$, we have $\|q-v(q)\| >
  \CR (Q)$.  Consider then a ball $\ov{B}_2(v(q),\|q-v(q)\|)$ centered at
  the vertex $v(q)$, with radius $\|q-v(q)\|$.  By definition of the
  circumcenter, $\CC(Q)$ must lie in the interior of
  $\ov{B}_2(v(q),\|q-v(q)\|)$. Consequently, $\| \CC (Q)-v (q) \| < \| q -
  v (q)\|$. Then, from $ \| \CC (Q) - v (q) \|^2 = \| \CC (Q) - q \|^2 + \|
  q - v (q) \|^2 - 2 (q-\CC(Q)) \cdot (q-v(q))$, we deduce
  \begin{equation*}
    2 (q-\CC(Q)) \cdot (q-v(q)) - \| \CC (Q) - q \|^2 
    = \| q - v (q) \|^2 - \| \CC (Q) - v (q) \|^2 > 0 \, ,
  \end{equation*}
  which implies the desired result.  
  
  Let us now prove (iii).  If $q \in \IC (Q)$, the claim is obviously
  satisfied since $\LN[\sm_Q](q)=0$.  Assume then that $q \not \in \IC
  (Q)$.  By Proposition~\ref{prop:remarkable}, $0 \not \in \partial \sm_Q
  (q)$, and hence $\LN[\sm_Q](q) \neq 0$.  Let us prove $\LN[\sm_Q](q)
  \cdot n_e > 0$ reasoning by contradiction. If $\LN[\sm_Q](q) \cdot n_e
  \le 0$, then $d/dt \left( \D_e (q + t \LN[\sm_Q](q)) \right)_{t=0} = n_e
  \cdot \LN[\sm_Q](q) \le 0$, which implies that $\D_e (q + t
  \LN[\sm_Q](q)) \le \D_e (q) = \sm_Q (q)$ for $t>0$ small enough. On the
  other hand, invoking Theorem~\ref{th:descent-direction} for the function
  $-\sm_Q$, we have that $\sm_Q (q) + t \| \LN[\sm_Q](q)\|^2/2 \le \sm_Q (q
  + t \LN[\sm_Q](q)) \le \D_e (q + t \LN[\sm_Q](q))$.  Gathering both
  facts, we conclude $t \| \LN[\sm_Q](q)\|^2/2 \le 0$, which is a
  contradiction.
  
  Let us now prove (iv).  By definition, $\D_e (q) \le \IR (Q)$.  This
  inequality is strict if $q \not \in \IC(Q)$. Let $x \in \IC(Q)$. If we
  take a point $O$ in the edge $e$, then the function $\D_e$ can be
  expressed as $\D_e (p) = (p-O) \cdot n_e$.  Then, we have
  \begin{gather*}
    \D_e(x) = (x-O) \cdot n_e = (x-q) \cdot n_e + (q-O) \cdot n_e = (x-q)
    \cdot n_e + \D_e (q) \, .
  \end{gather*}
  Since $\D_e(x) \ge \sm_Q (x) =\IR (Q)$, we conclude that $(x-q) \cdot n_e
  \ge \IR(Q) - \D_e(q) \ge 0$, and that the second inequality is strict if
  $q \not \in \IC(Q)$.  \quad
\end{proof}

\section{Analysis of the multi-center functions}\label{sec:multi-center-analysis}

Here we study the locational optimization functions $\HHDC$ and $\HHSP$ for
the disk-covering and sphere-packing problems.  We characterize their
smoothness properties, generalized gradients, and critical points for
arbitrary numbers of generators.  

\subsection{Smoothness and generalized gradients}
We start by providing some alternative
expressions and useful quantities. We write
\begin{align*}
  \HHDC (P) = \max_{i\in\until{n}} G_i(P) \, , \quad \HHSP (P) =
  \min_{i\in\until{n}} F_i(P) \, ,
\end{align*}
where
\begin{align*}
  G_i(P) = \max_{q \in V_i(P)} \|q-p_i\| \, , \quad F_i(P) = \min_{ q \not \in
    \interior{V_i(P)}} \|q-p_i\| \, .
\end{align*}
%
% \sindex{G}{$G_i(P)$}{Largest distance from $p_i$ to the boundary of $V_i (P)$}
% \sindex{F}{$F_i(P)$}{Smallest distance from $p_i$ to the boundary of $V_i (P)$}
%
Note that $G_i (P) = \lg_{V_i(P)}(p_i)$ and $F_i(P)=\sm_{V_i(P)}(p_i)$,
where, for $i \in \until{n}$,
\begin{align*}
  \lg_{V_i}:V_i \rightarrow \real \, , \quad \sm_{V_i}: V_i \rightarrow
  \real \,.
\end{align*}
Proposition~\ref{prop:remarkable} provides an explicit expression for the
generalized gradients of $\lg_{V_i}$ and $\sm_{V_i}$ when the Voronoi cell
$V_i$ is held fixed. Despite the slight abuse of notation, it is convenient
to let $\partial \lg_{V_i(P)}(p_i)$ denote $\partial \lg_{V}(p_i)_{|
  V=V_i(P)}$, and let $\partial \sm_{V_i(P)}(p_i)$ denote $\partial
\sm_{V}(p_i)_{| V=V_i(P)}$.

In contrast to this analysis at fixed Voronoi partition, the properties of
the functions $G_i$ and $F_i$ are strongly affected by the dependence on
the Voronoi partition $\VV(P)$.  We endeavor to characterize these
properties in order to study $\HHDC$ and $\HHSP$.
\begin{proposition}\label{prop:locally-Lipschitz}
  The functions $G_i, F_i : Q^n \rightarrow \real$ are locally Lipschitz
  and regular. As a consequence, the locational optimization functions
  $\HHDC, \HHSP :Q^n \rightarrow \real$ are locally Lipschitz and regular.
\end{proposition}

\begin{prof}{\it Proof. (a) $G_i$ is locally Lipschitz and regular}.
  The definition of the function $G_i$ admits the following alternative
  expression
  \begin{align}
    \label{eq:alternative-G_i}
    G_i(P) & = \max_{v \in \VE(V_i)} \| p_i - v \| \, .
  \end{align}
  Let $P_0$ be nondegenerate at the $i$th generator.  Then, there exists a
  neighborhood $U$ of $P_0$ where the set $\NN (i)$ does not change.  Let
  $\{ v_1,\dots, v_{M_1} \}$, $\{ w_1,\dots, w_{M_2} \}$, $\{ z_1,\dots,
  z_{M_3} \}$ be the vertexes of $V_i$ of types (a), (b) and (c)
  respectively.  Then, $G_i$ can be locally written as
  \begin{align*}
    G_i(P) =  \max \left\{ 
      \max_{\ell \in \until{M_1}} \|v_\ell-p_i\|, 
      \max_{\ell \in \until{M_2}} \|w_\ell-p_i\|,  
      \max_{\ell \in \until{M_3}} \|z_\ell-p_i\| \right\}  ,
  \end{align*}
  for all $P\in U$. Therefore, $G_i$ restricted to $U$ coincides with the
  function $\map{\GG_{\NN(i)}}{Q^n}{\real}$ defined by
  \begin{equation}
    \label{eq:GG}
    \GG_{\NN(i)} (P) = \max \left\{ 
      \max_{\ell \in \until{M_1}} \|v_\ell-p_i\|, 
      \max_{\ell \in \until{M_2}} \|w_\ell-p_i\|,  
      \max_{\ell \in \until{M_3}} \|z_\ell-p_i\| \right\} .
  \end{equation}
  The function $\GG_{\NN(i)}$ is the maximum of a fixed finite set of
  locally Lipschitz and regular functions, and consequently, locally
  Lipschitz and regular by Proposition~\ref{prop:gradient-F-nice}. We
  conclude that $G_i$ is both locally Lipschitz and regular at $P_0$.
   
  Let $P_0$ be degenerate at the $i$th generator. Then, in any neighborhood
  $U$ of $P_0$ there are different sets of neighbors of the $i$th
  generator.  Indeed, because the number of generators, edges of the
  boundary $Q$ and vertexes of $Q$ is finite, there is only a finite number
  of different sets of neighbors of the $i$th generator over $U$, say
  $\NN^1 (i), \dots, \NN^L (i)$.  This implies that $G_i$ admits the
  following alternative expression over $U$
  \begin{align}\label{eq:GG-multiple}
    G_i (P) = \min \left\{ \GG_{\NN^1(i)} (P), \dots, \GG_{\NN^L(i)}(P)
    \right\} \, .
  \end{align}
  Again resorting to Proposition~\ref{prop:gradient-F-nice}, we conclude
  that $G_i$ is both locally Lipschitz and regular at $P_0$.
  
  \noindent  \emph{(b)  $F_i$  is locally   Lipschitz  and  regular}.
  From the definition of $F_i$, it is clear that its value at a
  configuration $P$ is attained at the boundary of the Voronoi region
  $V_i$.  Therefore, one only minimizes among the edges associated with the
  Voronoi neighbors $\NN (i)$ and the edges of $Q$ with non-empty
  intersection with $V_i$.  Moreover, one can also see that the minimum
  must be attained at a point of the form $\proj{e}(p_i)$, for some edge
  $e$ of $V_i$.  Now, consider the function
  $\map{\FF_{\NN(i)}}{Q^n}{\real}$ defined by
  \begin{equation} 
    \label{eq:expression-F_i}
    \FF_{\NN(i)}(P) =
    \min \left\{ \min_{j \in \until{n}} \| p_i - \frac{p_i+p_j}{2} \| \, , \;
      \min_{e \in \ED (Q)} \D_e(p_i) \right\}
  \end{equation}
  We shall prove that $\FF_{\NN(i)}$ coincides with $F_i$.  If $k \not \in
  \NN (i)$, then $(p_i+p_k)/{2} \not \in V_i$.  Since $Q \setminus V_i$ is
  open, there exists a neighborhood of $(p_i+p_k)/{2}$ such that $U \subset
  Q/V_i$. Therefore,
  \begin{equation*}
    \| p_i - \frac{p_i+p_j}{2} \| > \min_{q \in U}  \| p_i - q \| \ge
    \min_{q \not \in \interior{V_i}}  \| p_i - q \| = F_i (P) \, .  
  \end{equation*}
  If an edge $e$ of $Q$ does not intersect $V_i$, then
  $\proj{e} (p_i) \not \in V_i$. Using again the fact that
  $Q \setminus V_i$ is open, there exists a neighborhood $U$ of
  $\proj{e} (p_i)$ in $\real^2$ such that $U \cap Q \subset
  Q \setminus V_i$. Then,
  \begin{equation*}
    \| p_i - \proj{e} (p_i) \| >  \min_{q \in U \cap Q} \| p_i -
    q \| \ge F_i (P) \, .
  \end{equation*}
  As a consequence of the previous inequalities, $F_i$ equals
  $\FF_{\NN(i)}$. Being the minimum of a fixed finite number of locally
  Lipschitz and regular functions on $Q^n$, $F_i$ is also locally Lipschitz
  and regular by Proposition~\ref{prop:gradient-F-nice}.  \quad
\end{prof}

Next, one can actually prove the following stronger result.
\begin{proposition}\label{prop:globally-Lipschitz}
  The locational optimization functions $\HHDC, \HHSP :Q^n \rightarrow
  \real$ are globally Lipschitz, with Lipschitz constant equal to $1$.
\end{proposition}

\begin{prof}{\it Proof. (a) $\HHDC$ is globally Lipschitz}. 
  Let $P$, $P'$ be two configurations of the $n$ generators.  Without loss of
  generality, assume that $\HHDC(P) \le \HHDC(P')$.  Let $i$, $j$ and $q_0$,
  $q_0' \in Q$ be such that $\HHDC(P)=G_i(P) = \|q_0 - p_i \|$ and
  $\HHDC(P')=G_j(P') = \| q'_0 - p'_j \|$. Now, consider the set
  $B_2(q'_0,G_i(P))$.  Then there exists a $k$ such that $p_k \in
  \ov{B}_2(q'_0,G_i(P))$ (otherwise, $\| q'_0 - p_l \| > G_i(P)$, which
  contradicts the definition of the function $\HHDC$).  On the other hand, we
  necessarily have that $p'_k \not \in B_2(q'_0,G_j(P'))$, since otherwise
  $\| q'_0 - p'_k \| < \| q'_0 - p'_j \|$, which implies that $q'_0 \not \in
  V'_j$, contradiction.  Finally, we apply the triangle inequality to obtain
  $\| q'_0 - p'_k \| \le \| q'_0 - p_k \| + \| p_k - p'_k \|$. Gathering the
  previous facts, we have
  \begin{multline*}
    | \HHDC(P') - \HHDC(P)| = G_j(P') - G_i(P) \\
    \le \| q'_0 - p'_k \| - \| q'_0 - p_k \| \le \| p_k - p'_k \| \le
    \|P - P' \| \, .
  \end{multline*}

  \noindent \emph{(b) $\HHSP$ is globally Lipschitz}. Let $P$, $P'$ be two
  configurations of the $n$ generators.  Without loss of generality, assume
  that $\HHSP(P) \le \HHSP(P')$.  Let $i$, $j$ and $q_0$, $q_0' \in Q$ be
  such that $\HHSP(P)=F_i(P) = \| q_0 - p_i \|$ and $\HHSP(P') = F_j(P') =
  \| q'_0 - p'_j \|$.  We treat separately the following two cases: (i)
  $q_0$ does not belong to the boundary of $Q$, and (ii) $q_0$ belongs to
  the boundary of $Q$.  In case (i), it necessarily exists $k \in \NN (i)$
  such that $\| q_0-p_i \| = \| q_0-p_k \|$.  If $\| q_0-p'_i \| \ge F_j
  (P')$, then
  \begin{align}\label{eq:intermediate}
    |\HHSP(P') - \HHSP(P) |= F_j (P') - F_i (P) & \le \| q_0 - p_i' \| - \|
    q_0 - p_i \| \nonumber \\
    & \le \| p_i - p_i' \| \le \| P - P' \| \, .
  \end{align}
  If, on the contrary, $\| q_0-p'_i \| < F_j (P')$, then $q_0 \in
  \interior{V_i'}$. Therefore, $\| q_0 - p_k'\| \ge F_k (P') \ge F_j (P')$. Now,
  we perform the same computation as in~\eqref{eq:intermediate} to conclude
  $|\HHSP(P') - \HHSP(P) | \le \| P - P' \|$.
  
  In case (ii), we prove that $\| q_0 - p_i' \| \ge F_j (P')$. Suppose this
  is not true, i.e., $\| q_0 - p_i' \| < F_j (P')$. Let $m= q_0 + \eps (q_0 -
  p_i')$, with sufficiently small $\eps >0$ such that $\| m - p_i' \| < F_j
  (P')$. Clearly $m \not \in Q$. On the other hand, by definition
  $B_2(p_i',F_i(P')) \subset V_i'$. Now, we have,
  \[
  B_2(p_i',F_j(P')) \subset B_2(p_i',F_i(P')) \subset V_i' \subset Q \, .
  \] 
  But, since $\| m - p_i' \| < F_j (P')$, then $m \in B_2(p_i',F_j(P'))
  \subset Q$, which is a contradiction. Therefore, $\| q_0 - p_i' \| \ge F_j
  (P')$, and now the same argument as in~\eqref{eq:intermediate} guarantees
  that $|\HHSP(P') - \HHSP(P) | \le \| P - P' \|$.  \quad
\end{prof}

\begin{figure}[htb]
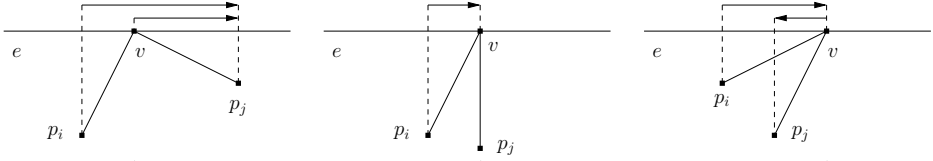
 
  \begin{center}
    \resizebox{.3\linewidth}{!}{\input{lambda-pos.tex}}\quad
    \resizebox{.3\linewidth}{!}{\input{lambda-zero.tex}}\quad
    \resizebox{.3\linewidth}{!}{\input{lambda-neg.tex}}
  \end{center}
  \caption{To illustrate eq.~\eqref{eq:def-lambda} we draw 
    the vectors $\proj{e}(p_j-v(e,i,j))$ and $\proj{e} (p_j-p_i)$ for
    various locations of $p_i$, $p_j$, and $e$. The left, center and right
    figures correspond to $\lambda(e,i,j)>0$, $\lambda(e,i,j)=0$,
    $\lambda(e,i,j)<0$, respectively. }
  \label{fig:lambda}
\end{figure}

We now introduce some quantities that are useful in characterizing the
generalized gradient of the functions $G_i$.  Given a vertex of type (b),
$v=v(e,i,j)$, determined by the edge $e$ and two generators $p_i$ and $p_j$,
we consider the scalar function $\lambda(e,i,j)$ defined by
\begin{equation} \label{eq:def-lambda}
  \proj{e} (p_j-v(e,i,j)) = \lambda(e,i,j) \, \proj{e} (p_j-p_i)
\end{equation}
where $P_e$ is the orthogonal projection onto the edge $e$; see
Fig.~\ref{fig:lambda}.  One can see that $\lambda(e,i,j)+\lambda(e,j,i)=1$.
%%
% \sindex{lambda}{$\lambda(e,i,j)$}{Scalar function associated with the vertex
%   $v(e,i,j)$}
%%
If $e$ is a segment in the line $ax+by+c=0$, $(\Delta x_{ij},\Delta y_{ij}) =
p_j - p_i$, $(x_m,y_m) =(p_i+p_j)/2$, then one can show
\begin{equation*}
  \lambda(e,i,j) = \frac{1}{2} - \frac{(a \Delta x_{ij} + b \Delta y_{ij})(a x_m + b y_m + c)}{(a
    \Delta y_{ij} - b \Delta x_{ij})^2} \, .
\end{equation*}
Given a vertex of type (a), $v=v(i,j,k)$, determined by the three generators
$p_i$, $p_j$, and $p_k$, we consider the scalar function $\mu(i,j,k)$ defined
by
\begin{equation*}
  \proj{e_{jk}} (p_\ell-v(i,j,k)) = \mu(i,j,k) \, \proj{e_{jk}} (p_\ell-p_i)
\end{equation*}
where $e_{jk}$ is the bisector of $p_j$ and $p_k$ and where $p_\ell=p_j$ if
$p_j$ belongs to the half-plane defined by $e_{jk}$ containing $p_i$, and
$p_\ell=p_k$ otherwise. One can see that $\mu(i,j,k)=\mu(i,k,j)$ and that
$\mu(i,j,k)+\mu(j,k,i)+\mu(k,i,j)=1$.
%%
% \sindex{mu}{$\mu(i,j,k)$}{Scalar function associated with the vertex
%   $v(i,j,k)$}
%%
From the expression for $\lambda$, one can obtain
\begin{equation*}
  \mu(i,j,k) = \frac{1}{2} + 
  \frac{
    (\Delta  x_{ij} \Delta x_{jk} + \Delta y_{ij} \Delta y_{jk})
    (\Delta x_{ik} \Delta x_{jk} + \Delta y_{ik} \Delta y_{jk})
    }%
  {2 (x_k \Delta y_{ij} - x_j \Delta y_{ik} + x_i \Delta y_{jk})^2} .
\end{equation*}
Note that, in general, $\lambda$ and $\mu$ are not positive functions.  Now
we are ready to describe in detail the structure of the generalized
gradient of the functions $G_i$, $F_i$.

\begin{proposition} \label{prop:gradient-G}
  The generalized gradient of $\map{G_i}{Q^n}{\real}$ at $P\in Q^n$ is
  \begin{align*}
    \partial G_i (P) & = \co \setdef{\partial_v G_i (P) \in (\real^2)^n}%
    { v\in\VE(V_i(P)) \; \text{such that} \; G_i (P) = \| p_i - v \|}
  \end{align*}
  where we consider separately the following cases. If $v=v(i,j,k)$ is a
  nondegenerate vertex of type (a), then
  \begin{multline*}
    \partial_{v(i,j,k)} G_i (P) =
    \partial_{v(k,i,j)} G_k (P) =
    \partial_{v(j,k,i)} G_j (P) = \\
     (0,\dots,
    \underbrace{\mu(i,j,k) \versus{p_i-v}}_{\text{$i$th place}}
    ,\dots,
    \underbrace{\mu(j,k,i) \versus{p_j-v}}_{\text{$j$th place}}
    ,\dots,
    \underbrace{\mu(k,i,j) \versus{p_k-v}}_{\text{$k$th place}}
    ,\dots, 0)
  \end{multline*}
  where, without loss of generality, we let $i<j<k$.  If $v=v(e,i,j)$ is a
  nondegenerate vertex of type (b), then
  \begin{multline*}
    \partial_{v(e,i,j)} G_i (P)=
    \partial_{v(e,j,i)} G_j (P) \\
    = (0,\dots,
    \underbrace{\lambda(e,i,j) \versus{p_i-v}}_{\text{$i$th place}}
    ,\dots,
    \underbrace{\lambda(e,j,i) \versus{p_j-v}}_{\text{$j$th place}}    
    ,\dots,0)
  \end{multline*}
  where, without loss of generality, we let $i<j$. If $v=v(e,f,i)$ is a
  nondegenerate vertex of type (c), then
  \begin{equation*}
    \partial_{v(e,f,i)} G_i (P) =
    (0,\dots,0,\underbrace{\versus{p_i-v}}_{\text{$i$th place}},0,\dots,0).
  \end{equation*}
  Finally, if the vertex $v$ is degenerate, i.e., if~$v$ is determined by
  $d>3$ elements (generators or edges), then there are $\binom{d-1}{2}$
  pairs of elements which determine the vertex~$v$ together with the
  generator~$p_i$. In this case, $\partial_v G_i (P)$ is the convex hull of
  $\partial_{v(\alpha,\beta,\gamma)} G_i (P)$ for all $\binom{d-1}{2}$ such
  triplets $(\alpha,\beta,\gamma)$.
\end{proposition}

Note that, at all nondegenerate configurations $P$, the quantity
$\partial_v G_i(P)$ is the generalized gradient of the function
$(p_1,\dots,p_n) \mapsto \| p_i-v(i,j,k)\|$; however, this interpretation
cannot be given when $P$ is degenerate.

\begin{proof}
  We present the proof for the expression for $\partial G_i (P)$. Let us
  consider first the case when $P$ is nondegenerate configuration for the
  $i$th generator.  According to the proof of
  Proposition~\ref{prop:locally-Lipschitz}, $G_i$ coincides with the
  function $\GG_{\NN(i)}$ over a neighborhood $U$ of $P$.  Hence, $\partial
  G_i (P) = \partial \GG_{\NN(i)} (P)$ which, according to
  eq.~\eqref{eq:GG} and Proposition~\ref{prop:gradient-F-nice}, takes
  the form
  \begin{equation*}
    \co \setdef{ \pder{}{P}  \| v-p_i \|}%
    { v \in \VE(V_i(P)) \;\text{such that}\; \| v-p_i \| = G_i (P) }
    \, .
  \end{equation*}  
  If $v=v(i,j,k)$ is a nondegenerate vertex of type (a), then
  \begin{align*}
    \pder{}{p_i} \| p_i - v(i,j,k) \| &= \versus{p_i - v} \left( I_2
      - \pder{v}{p_i} \right) = \mu (i,j,k) \versus{p_i - v} \, ,
    \\
    \pder{}{p_j} \| p_i - v(i,j,k) \| &= -\versus{p_i - v} \left(
      \pder{v}{p_j} \right) = \mu(j,k,i) \versus{p_j - v} \, ,
    \\
    \pder{}{p_\ell} \| p_i - v(i,j,k) \| &= 0 \, , \qquad \ell \not = i,j,k  \,.
  \end{align*}
  If $v=v(e,i,j)$ is a nondegenerate vertex of type (b), then
  \begin{align*}
    \pder{}{p_i} \| p_i - v(e,i,j) \| &= \versus{p_i - v} \left( I_2 -
      \pder{v}{p_i} \right) = \lambda(e,i,j) \versus{p_i - v} \, ,
    \\
    \pder{}{p_j} \| p_i - v(e,i,j) \| &= - \versus{p_i - v} \left(
      \pder{v}{p_j} \right) = \lambda(e,j,i) \versus{p_j - v} \, , 
    \\
    \pder{}{p_\ell} \| p_i - v(e,i,j) \| &= 0 \, , \qquad \ell \not = i,j \, .
  \end{align*}
  If $v=v(e,f,i)$ is a nondegenerate vertex of type (c), then
  \begin{align*}
     \pder{}{p_i} \| p_i - v(e,f,i) \| &= \versus{p_i - v} \, ,
     \\
     \pder{}{p_\ell} \| p_i - v(e,f,i) \| &= 0 \, , \qquad
     \ell \not= i \,.
  \end{align*}
 
  If $P$ is a degenerate configuration at the $i$th generator, then this
  function can be expressed as in eq.~\eqref{eq:GG-multiple} in a
  sufficiently small neighborhood $U$ of $P$. According to
  Proposition~\ref{prop:gradient-F-nice}, the generalized gradient of $G_i$
  is given by the convex hull of the generalized gradients of each of the
  functions $\GG_{\NN^1(i)}, \dots, \GG_{\NN^L(i)}$.  The claim now follows
  by reproducing the previous discussion for the generalized gradients of
  each of the functions $\GG_{\NN^\ell(i)}$, $\ell \in \until{L}$. \quad
\end{proof}

The expression for $\partial F_i (P)$ can be deduced in an analogous (and
simpler) way, since according to the proof of
Proposition~\ref{prop:locally-Lipschitz}, it is not necessary to establish
any distinction between the degenerate and the nondegenerate
configurations. Accordingly, we state the following result without proof.

\begin{proposition} \label{prop:gradient-F}
  The generalized gradient of $\map{F_i}{Q^n}{\real}$ at $P\in Q^n$ is
  \begin{align*}
    \partial F_i (P) & = \co \setdef{\partial_e F_i (P) \in (\real^2)^n}%
    {e\in\ED(V_i(P)) \; \text{such that} \; F_i (P) = \D_e(p_i)}
  \end{align*}
  where, if $e=e(i,j)$ is an edge of type (a), then
  \begin{equation*}
    \partial_{e(i,j)} F_i (P) =    \partial_{e(j,i)} F_j (P) =
    \frac{1}{2} (0,\dots,    
    \underbrace{n_{e(i,j)}}_{\text{$i$th place}}, \dots,
    \underbrace{-n_{e(i,j)}}_{\text{$j$th place}}, \dots,0),
  \end{equation*}
  and if $e=e(i)$ is an edge of type (b), then
  \begin{equation*}
    \partial_{e(i)} F_i (P)     =
    (0,\dots,\underbrace{n_{e(i)}}_{\text{$i$th place}},\dots,0) .
  \end{equation*}
\end{proposition}

Next, we give conditions under which the functions $\lambda$ and $\mu$ take
positive values.

\begin{lemma}\label{le:positive-eigenvalues}
  Let $P \in Q^n$ and let $v \in \VE_{\DC}(\VV(P))$.  Then,
  \begin{enumerate}
  \item if $v$ belongs to an edge $e$ of $Q$, then there exist generators
    $p_i$ and $p_j$ such that $\lambda (e,i,j)$ and $\lambda(e,j,i)$ are
    positive, and
  \item if $v$ belongs to $\interior{Q}$, then there exist generators
    $p_i$, $p_j$ and $p_k$ such that $\mu (i,j,k)$, $\mu (j,k,i)$ and $\mu
    (k,i,j)$ are positive.
  \end{enumerate}
\end{lemma}

\begin{proof}
  Consider first the case when $v$ is nondegenerate. If $v$ is in the edge
  $e$ of $Q$ (i.e. $v$ is of type (b)), let $p_i$ and $p_j$ the two
  generators determining it. From the definition of $\lambda$, one sees
  that the values $\lambda(e,i,j) =0$ and $\lambda(e,j,i)=0$ correspond to,
  respectively, $p_j$ and $p_i$ lying on the orthogonal line to $e$ passing
  through $v(e,i,j)$. If $\lambda(e,i,j) \le 0$, then there exists $w \in e
  \cap V_j$ such that $\| p_j - w\| > \| p_j - v\| = \HHDC (P)$, which is a
  contradiction. Therefore $\lambda(e,i,j) > 0$. The same argument
  guarantees $\lambda(e,j,i)>0$.  If $v$ is of type (a), and $p_i$, $p_j$
  and $p_k$ are the elements determining it, a similar argument leads to
  the conclusion that $\mu (i,j,k)$, $\mu (j,k,i)$ and $\mu (k,i,j)$ are
  positive.
  
  Consider the case when $v$ is degenerate.  Let $\{i_1,\dots,i_m \}$ be such
  that $v \in V_{i_j}$, $j \in \until{m}$. Assume $v$ is an edge $e$ of $Q$.
  Let $l$ denote the orthogonal line to the edge $e$ passing through $v$. We
  claim that there must exist generators in $\{i_1,\dots,i_m \}$ on both
  sides of $l$. Assume this is not the case, i.e., $\{p_{i_1},\dots,p_{i_m}
  \}$ are contained in one of the closed half-planes defined by $l$, say
  $l_-$.  Take $w \in l_+ \cap e$ arbitrarily close to $v$. Since
  $\{p_{i_1},\dots,p_{i_m} \} \subset l_-$, we have $\| p_{i_j} - w \| > \|
  p_{i_j} - v \|$. On the other hand, since no generator outside the set
  $\{p_{i_1},\dots,p_{i_m} \}$ is involved in the definition of $v$, there
  must exist $j^*$ such that $w \in V_{i_{j^*}}$.  Therefore, $G_{i_{j^*}}
  (P) \ge \| p_{i_{j^*}} - w \| > \| p_{i_{j^*}} - v \| = \HHDC (P)$, which
  is a contradiction. Assume now that $v \in \interior{Q}$. Our claim is
  that, for any line $l$ passing through $v$, there must exist generators on
  both sides of $l$. If this is not the case, i.e., $\{p_{i_1},\dots,p_{i_m}
  \} \subset l_-$, then take $w \in B_2(v,\eps) \cap l_+ \cap o$, where $o$
  denotes the orthogonal line to $l$ passing through $v$. As before, $w \in
  V_{i_{j^*}}$ for some $j^*$ and $\| p_{i_{j^*}} - w \| > \| p_{i_{j^*}} - v
  \|$, which yields a contradiction.  \quad
\end{proof}

This completes our analysis of the generalized gradients of $G_i$ and $F_i$
and, with these results, we return to studying the generalized gradients of
$\HHDC$ and $\HHSP$.  An immediate consequence of
Propositions~\ref{prop:gradient-F-nice} and~\ref{prop:locally-Lipschitz} is
that
\begin{align}\label{eq:generalized-gradients}
  \partial \HHDC (P) &= \co \setdef{\partial G_i (P)}{i \in I(P)} \, ,
  \nonumber \\
  \partial \HHSP (P) &= \co \setdef{\partial F_i (P)}{i \in I(P)} \, .
\end{align}
Furthermore we can provide the following more detailed characterization.

\begin{proposition}\label{prop:proj-generalized-gradients}
  Let $P \in Q^n$.  For each $i \in \until{n}$, the image by $\pi_i$ of the
  generalized gradients of $\HHDC$ and $\HHSP$ at $P$ is given by
  \begin{align*}
    \pi_i (\partial \HHDC (P)) & =
      \begin{cases}
        \pi_i (\partial G_i (P)) & \; \hbox{if} \; \; i \in I(P) \, , \;
        \VE_{\DC}(\VV (P)) \subset \VE (V_i(P)) \\
        \co \{ \pi_i (\partial G_i (P)), 0 \} & \; \hbox{if} \; \; i \in I(P)
        \, , \; \VE_{\DC}(\VV (P)) \not \subset \VE (V_i(P)) \\
        0 & \; \hbox{if} \; \; i \not \in I(P) \\
      \end{cases}
      \\
      \pi_i (\partial \HHSP (P)) & =
      \begin{cases}
        \pi_i (\partial F_i (P)) & \; \hbox{if} \; \; i \in I(P) \, , \;
        \ED_{\SP}(\VV (P)) \subset \ED (V_i(P)) \\
        \co \{ \pi_i (\partial F_i (P)),0 \} & \; \hbox{if} \; \; i \in I(P)
        \, , \;  \ED_{\SP}(\VV (P)) \not \subset \ED (V_i(P)) \\
        0 & \; \hbox{if} \; \; i \not \in I(P) \\
      \end{cases}
  \end{align*}
\end{proposition}

\begin{proof}
  From eq.~\eqref{eq:generalized-gradients}, if $i \not \in I(P)$, then
  $\pi_i (\partial \HHDC (P)) = 0$, $\pi_i (\partial \HHSP (P))=0$.  If $i
  \in I(P)$, then using Proposition~\ref{prop:gradient-G}, we deduce that
  the generators $p_j$ such that $\partial G_j$ has a nonzero entry in the
  $i$th place (and hence contributes to the projection by $\pi_i$ of
  $\partial \HHDC$) must share a vertex with the $i$th generator.
  Analogously, if $i \in I(P)$, then using
  Proposition~\ref{prop:gradient-F}, we deduce that the generators $p_j$
  such that $\partial F_j$ has a nonzero entry in the $i$th place (and
  hence contributes to the projection by $\pi_i$ of $\partial \HHSP$) must
  satisfy $j \in \NN(i)$.  For the disk-covering function, if $v$ is a
  common vertex of $V_i$ and $V_j$, determined by $i$, $j$ and a third
  element $\alpha$, then $\partial_{v(\alpha,j,i)} G_j =
  \partial_{v(\alpha,i,j)} G_i$, and the expression for $\pi_i (\partial
  \HHDC (P))$ then follows. The argument for the expression of $\pi_i
  (\partial \HHSP (P))$ is analogous.  \quad
\end{proof}

\subsection{Critical points}
Having characterized the generalized gradients of $\HHDC$ and $\HHSP$, we
now turn to studying their critical points.
\begin{theorem}[Minima of $\HHDC$] \label{th:minima-DC}
  Let $P \in Q^n$ be a nondegenerate configuration and $0 \in \interior{\partial
    \HHDC(P)}$.  Then, $P$ is a strict local minimum of $\HHDC$, all
    generators are active and $P$ is a circumcenter Voronoi configuration.
\end{theorem}

\begin{proof}
  Since $P$ is nondegenerate, note from Proposition~\ref{prop:gradient-G} that
  $\partial_v G_i$ is a singleton for each $v \in \VE (V_i(P))$, $i \in
  \until{n}$.  Let $w \in (\real^2)^n$. We claim that moving the
  configuration of the generators from $P$ in the direction $w$ can only
  increase the cost.  The hypothesis $0 \in \interior{\partial \HHDC(P)}$
  implies by Lemma~\ref{le:positively-span-iff-I} that there exists $i$ and
  $v \in \VE(V_i (P)) \cap \VE_{\DC}(\VV (P))$ such that $w \cdot \partial_v
  G_{i}(P) >0$.  Since $P$ is nondegenerate, $v$ will still belong to $V_{i}
  (P + \eps w)$ for sufficiently small $\eps>0$, and consequently $ \HHDC (P
  + \eps w) \ge G_{i} (P + \eps w) > G_{i}(P) = \HHDC (P)$. Therefore $P$ is
  a strict local minimum.
  
  Since $\pi_i$ is an open map, the set $\pi_i (\interior{\partial \HHDC
    (P)})$ is open for each $i \in \until{n}$. Therefore, $\pi_i
  (\interior{\partial \HHDC (P )}) \not = 0$, and hence all generators are
  active, i.e. $I(P ) = \until{n}$.  Let us see that all generators must also
  be centered. Assume $P$ is nondegenerate and consider the $i$th generator.
  Take $w \in \real^2$ and let $\ov{w} \in (\real^2)^n$ be the vector with
  has $w$ in the $i$th place and $0$ otherwise.  By
  Lemma~\ref{le:positively-span-iff-I}, there exist $j$ and $v \in \VE (V_j
  (P)) \cap \VE_{\DC}(\VV (P))$ such that $\ov{w} \cdot \partial_v G_j > 0$.
  Since $\ov{w} \cdot \partial_v G_j = w \cdot \pi_i (\partial_v G_j) > 0$,
  then necessarily $\pi_i (\partial_v G_j) \neq 0$, and therefore $v \in
  V_i(P)$ and $\pi_i (\partial_v G_j) = \pi_i (\partial_v G_i)$. The vertex
  $v$ is determined by $p_i$, $p_j$ and a third element, say $\alpha$.
  Depending on whether $\alpha$ corresponds to an edge or to another
  generator, we have that $\pi_i (\partial_v G_i)$ is equal to
  $\lambda(\alpha,i,j) \versus{p_i - v}$ or $\mu (\alpha,i,j) \versus{p_i -
    v}$. In any case, from Lemma~\ref{le:positive-eigenvalues}, we deduce
  that $\lambda(\alpha,i,j)$ (respectively $\mu (\alpha,i,j)$) belongs to
  $(0,1)$.  Therefore $ w \cdot \pi_i (\partial_v G_i) >0$ implies $w \cdot
  \versus{p_i - v} >0$. Consequently, $0 \in \interior{\partial
    \lg_{V_i(P)}(p_i)}$. By Proposition~\ref{prop:remarkable}, this implies
  that $p_i = \CC(V_i)$.  Hence, $P$ is a circumcenter Voronoi configuration.
  \quad
\end{proof}

\begin{theorem}[Maxima of $\HHSP$]\label{th:maxima-SP}
  Let $P \in Q^n$ and $0 \in \interior{\partial \HHSP(P)}$.  Then, $P$ is a
  strict local maximum of $\HHSP$, all generators are active and $P$ is a
  generic incenter Voronoi configuration.
\end{theorem}

\begin{proof}
  The proof of this result is analogous to the proof of
  Theorem~\ref{th:minima-DC}. Note that $0 \in \interior{\partial
  \sm_{V_i(P)}(p_i)}$ implies, by Proposition~\ref{prop:remarkable}, that $\IC
  (V_i(P)) = \{ p_i\}$, and hence $P$ is a generic incenter Voronoi
  configuration.  \quad
\end{proof}

\begin{figure}
  \begin{center}
    \resizebox{.3\linewidth}{!}{\input{local-minimum.tex}} \qquad\qquad
    \resizebox{.45\linewidth}{!}{\input{local-maximum-2.tex}}
    \caption{Local extrema of the disk-covering and the sphere-packing
      functions in a convex polygonal environment. The configuration on the
      left corresponds to a local minimum of $\HHDC$ with $0 \in \partial
      \HHDC (P)$ and $\interior{\partial \HHDC (P)} = \emptyset$. The
      configuration on the right corresponds to a local maximum of $\HHSP$
      with $0 \in \partial \HHSP (P)$ and $\interior{\partial \HHSP (P)} =
      \emptyset$. In both configurations, the $4$th generator is inactive
      and non-centered.}
    \label{fig:local-min-max}
  \end{center}
\end{figure}

\begin{remark}{\rm
    Theorems~\ref{th:minima-DC} and~\ref{th:maxima-SP} precisely provide
    the interpretation of the multicenter problems that we gave in
    Section~\ref{se:locational-functions}: since all generators are active,
    they share the same radius. If one drops the hypothesis that $0$
    belongs to the generalized gradient of the locational optimization
    function, then one can think of simple examples where $P$ is a local
    minimum of $\HHDC$ (respectively  local maximum of $\HHSP$), and there are
    generators which are inactive and non-centered, see
    Fig.~\ref{fig:local-min-max}.}
\end{remark}

\section{Dynamical systems for the multi-center problems}
\label{sec:multi-center-design}

In this section, we describe three algorithms that (locally) extremize the
multi-center functions for the disk-covering and the sphere-packing problems.
We first examine the gradient flow descent associated with the locational
optimization functions $\HHDC$ and $\HHSP$.  This flow is guaranteed to find
a local critical point, but it has the drawback of being centralized, as we
describe later.  Then, we propose two decentralized flows for each problem.
One roughly consists of a distributed implementation of the gradient descent.
As we show, it is very much in the spirit of behavior-based robotics.  The
other one follows the logical strategy given the result in
Theorems~\ref{th:minima-DC} and~\ref{th:maxima-SP}: each generator moves
toward the circumcenter (alternatively, incenter set) of its own Voronoi
polygon.  We call them Lloyd flows, since they resemble the original Lloyd
algorithm for vector quantization problems, where each quantizer moves toward
the centroid or center of mass of its own Voronoi region,
see~\cite{QD-VF-MG:99,RMG-DLN:98,SPL:82}. We present continuous-time versions
of the algorithms and discuss their convergence properties. In our setting,
the generators' location obeys a first order dynamical behavior described by
\begin{equation}\label{eq:dynamical-system}
  \dot{p}_i = u_i (p_1,\dots,p_n) \, , \quad i \in \until{n} \, .
\end{equation}
The dynamical system~\eqref{eq:dynamical-system} is said to be (strongly)
\emph{centralized} if there exists at least an $i\in \until{n}$ such that
$u_i(p_1,\dots,p_n)$ cannot be written as a function of the form $u_i
(p_i,p_{i_1},\dots,p_{i_m})$, with $m < n-1$.  The dynamical
system~\eqref{eq:dynamical-system} is said to be \emph{Voronoi-distributed}
if each $u_i(p_1,\dots,p_n)$ can be written as a function of the form $u_i
(p_i,p_{i_1},\dots,p_{i_m})$, with $i_k \in \NN (P,i)$, $k \in \until{m}$.
Finally, the dynamical system~\eqref{eq:dynamical-system} is said to be
\emph{nearest-neighbor-distributed} if each $u_i(p_1,\dots,p_n)$ can be
written as a function of the form $u_i (p_i,p_{i_1},\dots,p_{i_m})$, with
$\| p_i - p_{i_k}\| \le \| p_i - p_{j}\|$ for all $j \in \until{n}$, and $k
\in \until{m}$. A nearest-neighbor-distributed dynamical system is also
Voronoi-distributed.

It is well known that there are at most $3n-6$ neighborhood relationships
in a planar Voronoi diagram~\cite[see Section 2.3]{AO-BB-KS-SNC:00}.
Therefore, the number of Voronoi neighbors of each site is on average less
than or equal to $6$.  (Recall that sites are Voronoi-neighbors if they
share an edge, not just a vertex.) We refer to~\cite{JC-SM-TK-FB:02j} for
more details on the distributed character of Voronoi neighborhood
relationships.

Note that the set of indexes $\{i_1,\dots,i_m\}$ for an specific generator
$p_i$ of a Voronoi-distributed or a nearest-neighbor-distributed dynamical
system is not the same for all possible configurations $P$. In other words,
the identity of both the Voronoi neighbors and the nearest neighbors might
change along the evolution, i.e., the topology of the dynamical system is
\emph{dynamic}.

\subsection{Nonsmooth gradient dynamical systems}

Consider the (signed) generalized gradient descent
flow~\eqref{eq:natural-gradient} for the locational optimization functions
$\HHDC$ and $\HHSP$,
\[
\dot{P} = - \LN (\partial\HHDC) (P) \, , \quad \dot{P} = \LN (\partial \HHSP)
(P) \, .
\]
Alternatively, we may write for each $i \in \until{n}$,
\begin{align}
  \dot{p}_i &= - \pi_i ( \LN (\partial \HHDC) (p_1,\dots,p_n)) \, ,
  \label{eq:gradient-flow-continuous-DC} \\
  \dot{p}_i &= \pi_i ( \LN(\partial \HHSP) (p_1,\dots,p_n)) \, .
  \label{eq:gradient-flow-continuous-SP}
\end{align}
As noted in Section~\ref{se:stability-via-Lipschitz}, these vector fields are
discontinuous, and therefore their solution must be understood in the
Filippov sense.  Eq.~\eqref{eq:generalized-gradients} and
Propositions~\ref{prop:gradient-G} and~\ref{prop:gradient-F} provide an
expression of the generalized gradients at $P$, $\partial \HHDC (P)$ and
$\partial \HHSP (P)$.  One needs to first compute the generalized gradient,
then compute the least-norm element, and finally project to each of the $n$
components; therefore the expressions in
Proposition~\ref{prop:proj-generalized-gradients} are not helpful.  Note that
the least-norm element of convex sets can be computed efficiently,
see~\cite{SB-LV:02}, however closed-form expressions are not available in
general.

One can see that the compact set $Q^n$ is strongly invariant for both vector
fields $-\LN(\partial\HHDC)$ and $\LN(\partial \HHSP)$.  Regarding $-
\LN(\partial \HHDC)$, this is a consequence of
Proposition~\ref{prop:gradient-G} and of Lemma~\ref{le:positive-eigenvalues}.
Regarding $ \LN(\partial \HHSP)$, this is a consequence of
Proposition~\ref{prop:gradient-F}.

\begin{proposition}
  For the dynamical system~\eqref{eq:gradient-flow-continuous-DC}
  (respectively~\eqref{eq:gradient-flow-continuous-SP}), the generators'
  location $P=(p_1,\dots,p_n)$ converges asymptotically to the set of
  critical points of $\HHDC$ (respectively, of $\HHSP$).
\end{proposition}

\begin{proof}
  From Propositions~\ref{prop:locally-Lipschitz}
  and~\ref{prop:globally-Lipschitz}, we know that $\HHDC$ and $\HHSP$ are
  globally Lipschitz and regular over $Q^n$. The result follows from
  Proposition~\ref{prop:natural-gradient} considering the dynamical system
  restricted to the strongly invariant and compact domain $Q^n$.  \quad
\end{proof}

\begin{remark} {\rm
    The gradient dynamical systems enjoy convergence guarantees, but their
    implementation is centralized because of two reasons. First, all
    functions $G_i(P)$ (respectively $F_i(P)$) need to be compared in order to
    determine which generator is active. Second, the least-norm element of
    the generalized gradients depends on the relative position of the
    active generators with respect to each other and to the environment. }
\end{remark}

\begin{remark} {\rm
    As illustrated in Fig.~\ref{fig:pcenter-2} the evolution of the
    gradient dynamical systems may not leave fixed even the generators that
    are centers (circumcenter or incenters).
  \begin{figure}[htb]
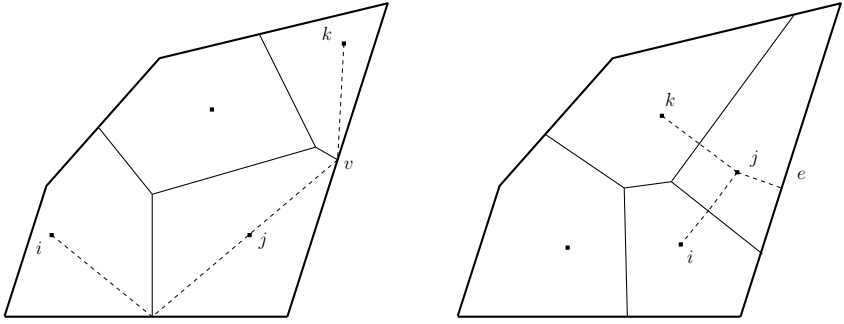

    \begin{center}
      \resizebox{.4\linewidth}{!}{\input{gradient-DC.tex}} \qquad
      \resizebox{.4\linewidth}{!}{\input{gradient-SP.tex}}
      \caption{Illustration of the gradient descent. In the left figure, the
        $j$th generator is in the circumcenter of its own Voronoi region, but
        the control law~\eqref{eq:gradient-flow-continuous-DC} drives it toward
        the vertex $v$. In the right figure, the $j$th generator is in the
        incenter of its own Voronoi region, but the control
        law~\eqref{eq:gradient-flow-continuous-SP} drives it away from the edge
        $e$.}\label{fig:pcenter-2}
      \end{center}
    \end{figure}  }
\end{remark}

\subsection{Nonsmooth dynamical systems based on distributed gradients}

In this section, we propose a distributed implementation of the previous
gradient dynamical systems and explore its relation with behavior-based rules
in multiple-vehicle coordination. Consider the following modifications of the
gradient dynamical
systems~\eqref{eq:gradient-flow-continuous-DC}-\eqref{eq:gradient-flow-continuous-SP},
\begin{align}
  \dot{p}_i &= - \LN (\partial \lg_{V_i(P)}) (P) \, ,
  \label{eq:distributed-gradient-flow-DC} \\ 
  \dot{p}_i &= \LN (\partial \sm_{V_i(P)}) (P) \, ,
  \label{eq:distributed-gradient-flow-SP}
\end{align}
for $i \in \until{n}$. Note that the
system~\eqref{eq:distributed-gradient-flow-DC} is Voronoi-distributed,
since $\LN (\partial \lg_{V_i(P)}) (P)$ is determined only by the
position of $p_i$ and of its Voronoi neighbors $\NN (P,i)$. On the
other hand, the system~\eqref{eq:distributed-gradient-flow-SP} is
nearest-neighbor-distributed, since $\LN ( \partial \sm_{V_i(P)})(P)$
is determined only by the position of $p_i$ and its nearest neighbors.

For future reference, let $\LN (\partial \lg_{\VV}) (P) = (\LN (\partial
\lg_{V_1 (P)}) (P),\dots,\LN (\partial \lg_{V_n(P)}) (P))$, $\LN (\partial
\sm_{\VV}) (P) = (\LN (\partial \sm_{V_1 (P)}) (P),\dots,\LN (\partial
\sm_{V_n (P)}) (P))$, and write
\[
\dot{P} = - \LN (\partial \lg_{\VV}) (P) \, , \quad \dot{P} = \LN ( \partial
\sm_{\VV}) (P) \, .
\]
As for the previous dynamical systems, note that these vector fields are
discontinuous, and therefore their solutions must be understood in the
Filippov sense. One can see that the compact set $Q^n$ is strongly invariant
for both vector fields $-\LN( \partial \lg_{\VV})$ and $\LN (\partial
\sm_{\VV})$.  This fact is a consequence of the expressions for the
generalized gradients of $\lg$ and $\sm$ in
Proposition~\ref{prop:remarkable}.  Note that in the 1-center
case,~\eqref{eq:gradient-flow-continuous-DC}
(respectively~\eqref{eq:gradient-flow-continuous-SP}) coincides
with~\eqref{eq:distributed-gradient-flow-DC} (respectively
with~\eqref{eq:distributed-gradient-flow-SP}).

\begin{proposition}
  Let $P \in Q^n$.  Then the solutions of the dynamical
  systems~\eqref{eq:distributed-gradient-flow-DC}
  and~\eqref{eq:distributed-gradient-flow-SP} starting at $P$ are
  unique.
\end{proposition}

\begin{prof}{\it Proof. (a) Uniqueness of solution for~\eqref{eq:distributed-gradient-flow-DC}}. 
  Let $D_{\lg}$ be the set of $P \in Q^n$ such that $P$ is nondegenerate and
  $\lg_{V_i(P)}(p_i)$ is attained at a single vertex for all $i$.  Note that
  $Q^n \setminus D_{\lg}$ has measure zero, and that the vector field $-\LN (
  \partial \lg_{\VV})$ is differentiable (and hence locally Lipschitz) when
  restricted to any connected component of~$D_{\lg}$. Let $P$, $P'$ belong to
  different connected components of $D_{\lg}$, and let $\|P- P'\| \le \eps$.
  Consider all the indexes $i$ at which the values of $\lg_{V_i(P)}(p_i)$ and
  $\lg_{V_i(P')}(p_i')$ are attained at different vertexes. For these
  indexes,
  \begin{align*}
    -\LN (\partial \lg_{V_i(P)}) (p_i) + \LN ( \partial\lg_{V_i(P')}) (p_i')
    = \versus{v-p_i} - \versus{w'-p_i'} \, ,
  \end{align*}
  for certain vertexes $v$ and $w'$. Note that for $\eps$ small enough, the
  vertex $w'$ in the Voronoi configuration $P'$ corresponds to a vertex $w$
  in the Voronoi configuration $P$.  By construction, $p_i$ and $p_i'$ belong
  to an $O(\eps)$ neighborhood of the bisector $b_{vw}$ determined by $v$ and
  $w$, and $n_{vw} \cdot (p_i -p_i') <0$. In addition, the component of $
  \versus{v-p_i} - \versus{w' - p_i'}$ along $b_{vw}$ is $O(\eps)$ whereas
  $n_{vw} \cdot \versus{v-p_i} >0$ and $n_{vw} \cdot \versus{w'-p_i'} =
  n_{vw} \cdot \versus{w-p_i} + O(\eps)$, with $ n_{vw} \cdot \versus{w-p_i}
  <0$. Then,
  \begin{multline*}
    \versus{v-p_i} - \versus{w'-p_i'}    \\
    = \proj{n_{vw}} (\versus{v-p_i} - \versus{w' - p_i'} +
    \proj{b_{vw}} (\versus{v-p_i} - \versus{w'-p_i'} \\
    = \proj{n_{vw}} (\versus{v-p_i} - \versus{w'-p_i'} + O (\eps) \, ,
  \end{multline*}
  and, in turn, for sufficiently small $\eps$
  \begin{multline*}
    (p_i - p_i') \cdot (\versus{v-p_i} - \versus{w'-p_i'}) \\
    = (n_{vw} \cdot (p_i -p_i')) (n_{vw} \cdot (\versus{v-p_i} -
    \versus{w'-p_i'})) + O(\eps^2) < 0 \, .
  \end{multline*}
  The result now follows from Theorem~1 at page~106 in~\cite{AFF:88}.

  \noindent \emph{(b) Uniqueness of solution
    for~\eqref{eq:distributed-gradient-flow-SP}}.  Let $D_{\sm}$ be the set
  of $P \in Q^n$ such that $\sm_{V_i(P)}(p_i)$ is attained at a single edge
  for all $i$. Note that $Q^n \setminus D_{\sm}$ has measure zero, and that
  the vector field $\LN ( \partial \sm_{\VV})$ is differentiable(and hence
  locally Lipschitz) when restricted to any connected component of~$D_{\sm}$.
  Let $P$, $P'$ belong to different connected components of $D_{\sm}$, and
  let $\|P- P'\| \le \eps$. Consider all the indexes $i$ at which the values
  of $\sm_{V_i(P)}(p_i)$ and $\sm_{V_i(P')}(p_i')$ are attained at different
  edges. Assume these edges are of type (a) (the type (b) case can be treated
  analogously). For these indexes,
  \begin{align*}
    \LN ( \partial \sm_{V_i(P)}) (p_i) - \LN ( \partial \sm_{V_i(P')}) (p_i')
    = \versus{p_i - p_j} - \versus{p_i' - p_k'} \, ,
  \end{align*}
  for some uniquely determined $p_j$ and $p_k'$, with $j \not = k$.  By
  construction, $p_i$ and $p_i'$ belong to an $O(\eps)$ neighborhood of the
  bisector $b_{jk}$ determined by $p_j$ and $p_k$, and $n_{kj} \cdot (p_i
  -p_i') <0$. In addition, the component of $ \versus{p_i - p_j} -
  \versus{p_i' - p_k'}$ along $b_{jk}$ is $O(\eps)$ whereas $n_{kj} \cdot
  \versus{p_i-p_j} >0$ and $n_{kj} \cdot \versus{p_i'-p_k'} = n_{kj} \cdot
  \versus{p_i-p_k} + O(\eps)$, with $ n_{kj} \cdot \versus{p_i-p_k} <0$.
  Then,
  \begin{multline*}
    \versus{p_i - p_j} - \versus{p_i' - p_k'} \\
    = \proj{n_{kj}} (\versus{p_i - p_j} - \versus{p_i' - p_k'} +
    \proj{b_{jk}}
    (\versus{p_i - p_j} - \versus{p_i' - p_k'} \\
    = \proj{n_{kj}} (\versus{p_i - p_j} - \versus{p_i - p_k} + O
    (\eps) \, ,
  \end{multline*}
  and, in turn, for sufficiently small $\eps$
  \begin{multline*}
    (p_i - p_i') \cdot (\versus{p_i - p_j} - \versus{p_i' - p_k'}) \\
    = (n_{kj} \cdot (p_i -p_i')) (n_{kj} \cdot (\versus{p_i - p_j} -
    \versus{p_i - p_k})) + O(\eps^2) < 0 \, .
  \end{multline*}
  The result now follows from Theorem~1 at page~106 in~\cite{AFF:88}. \quad
\end{prof}

\begin{remark}[Relation with behavior-based robotics: move toward the
  furthest-away vertex] {\rm The distributed gradient control law in
    the disk-covering setting~\eqref{eq:distributed-gradient-flow-DC}
    has an interesting interpretation in the context of behavior-based
    robotics.  Consider the $i$th generator.  If the maximum of
    $\lg_{V_i(P)}$ is attained at a single vertex $v$ of its Voronoi
    cell $V_i$, then $\lg_{V_i(P)}$ is differentiable at that
    configuration, and its derivative corresponds to $\versus{p_i-v}$.
    Therefore, the control law~\eqref{eq:distributed-gradient-flow-DC}
    corresponds to the behavior ``move toward the furthest vertex in
    own Voronoi cell.''  If there are two or more vertexes of $V_i$
    where the value $\lg_{V_i(P)}(p_i)$ is attained,
    then~\eqref{eq:distributed-gradient-flow-DC} provides an average
    behavior by computing the least-norm element in the convex hull of
    all $\versus{p_i-v}$ such that $\| p_i -v \| = \lg_{V_i(P)}(p_i)$.
  }
\end{remark}

\begin{remark}[Relation with behavior-based robotics: move away from the
  nearest neighbor] {\rm The distributed gradient control law in the
    sphere-packing setting~\eqref{eq:distributed-gradient-flow-SP} has also
    an interesting interpretation.  For the $i$th generator, if the minimum
    of $\sm_{V_i(P)}$ is attained at a single edge $e$, then $\sm_{V_i(P)}$
    is differentiable at that configuration, and its derivative is $n_e$.
    The control law~\eqref{eq:distributed-gradient-flow-SP} corresponds to
    the behavior ``move away from the nearest neighbor'' (where a neighbor
    can also be the boundary of the environment).  If there are two or more
    edges where the value $\sm_{V_i(P)}(p_i)$ is attained,
    then~\eqref{eq:distributed-gradient-flow-SP} provides an average
    behavior in an analogous manner as before. }
\end{remark}

\begin{proposition}\label{prop:distributed-gradient-DC}
  For the dynamical system~\eqref{eq:distributed-gradient-flow-DC}, the
  generators' location $P=(p_1,\dots,p_n)$ converges asymptotically to the
  largest weakly invariant set contained in the closure of $A_{\DC} (Q) =
  \setdef{ P \in Q^n}{i \in I(P) \implies p_i = \CC (V_i)}$.
%
%   \sindex{ADC}{$A_{\DC}(Q)$}{Set of configurations $P \in Q^n$ where all
%     active generators are in the circumcenter of its own Voronoi region}
%
\end{proposition}

\begin{proof}
  Let $a \in \setLieder{-\LN (\partial \lg_{\VV})}{\HHDC}(P)$. By definition,
  $a= -\LN ( \partial \lg_{\VV}) (P) \cdot \zeta$, for all $\zeta \in
  \partial \HHDC (P)$. Let $v \in \VE_{\DC}(\VV(P))$. From
  Lemmas~\ref{prop:gradient-G} and~\ref{le:positive-eigenvalues}, we know
  that, independently of the degenerate/nondegenerate character of the
  Voronoi partition at $v$, there always exist either an edge $e$ of $Q$ and
  generators $p_i$ and $p_j$, or generators $p_i$, $p_j$ and $p_k$, such that
  $\lambda (e,i,j)$, $\lambda (e,j,i) > 0$ (respectively $\mu (i,j,k)$,
  $\mu(j,k,i)$, $\mu (k,i,j)> 0$).  If $v$ is a vertex of type (b), then
  \begin{align}\label{eq:derivation}
    a &= -\LN (\partial \lg_{\VV}) (P) \cdot \partial_{v} G_{i}
    \\
    &= -\LN (\partial \lg_{V_i(P)}) (P) \cdot \lambda(e,i,j) \versus{p_i - v}
    -\LN (\partial \lg_{V_j (P)}) (P) \cdot \lambda(e,j,i) \versus{p_j - v} .
    \nonumber
  \end{align}
  From Lemma~\ref{le:simple-facts-DC-SP}(i) we conclude that $a \le 0$, and
  the inequality is strict if either $p_i \not = \CC (V_i)$ or $p_j \not =
  \CC (V_j)$.  The same conclusion can be derived if $v$ is a vertex of type
  (a).  Therefore, $\max \setLieder{-\LN (\partial \lg_{\VV})}{\HHDC}(P) \le
  0$ or $\setLieder{-\LN (\partial \lg_{\VV})}{\HHDC}(P) = \emptyset$. Now,
  resorting to the LaSalle principle (Theorem~\ref{th:LaSalle}), we deduce
  that the solution $P:[0,+\infty) \rightarrow Q^n$ starting from $P_0$
  converges to the largest weakly invariant set contained in $\ov{Z}_{-\LN (
    \partial \lg_{\VV}),\HHDC} \cap \HHDC^{-1} (\le \HHDC(P_0),P_0) \cap
  Q^n$.
  
  Let us see that $Z_{-\LN (\partial \lg_{\VV}),\HHDC}\cap Q^n$ is equal to
  $A_{\DC} (Q)$.  Take a configuration $P \in A_{\DC} (Q)$.  Then, $\LN
  (\partial \lg_{V_i(P)}) (P) =0$ if $i \in I(P)$, and $\pi_i (\zeta) =0$ if
  $i \not \in I(P)$, for any $\zeta \in \partial \HHDC (P)$ (cf.
  Proposition~\ref{prop:proj-generalized-gradients}). Consequently, $0 = -\LN
  (\partial \lg_{\VV}) (P) \cdot \zeta$, for all $\zeta \in \partial \HHDC
  (P)$, and so $0 \in \setLieder{-\LN (\partial \lg_{\VV})}{\HHDC}(P)$.
  Therefore, $A_{\DC}(Q) \subset Z_{-\LN (\partial \lg_{\VV}),\HHDC}$.  Now,
  consider $P \in Z_{-\LN (\partial \lg_{\VV}),\HHDC}$.  Then, $0 \in
  \setLieder{-\LN (\partial \lg_{\VV})}{\HHDC}(P)$, that is, $0 = -\LN
  (\partial \lg_{\VV}) (P) \cdot \zeta$, for all $\zeta \in \partial \HHDC
  (P)$. If $P$ is nondegenerate, we deduce from eq.~\eqref{eq:derivation} and
  Lemma~\ref{le:simple-facts-DC-SP} that all the active generators are
  centered, i.e., $P \in A_{\DC}(Q)$. If $P$ is degenerate, consider a
  degenerate vertex $v$ where the value of $\HHDC (P)$ is attained.  For
  simplicity, we deal with the case where $v$ is contained in an edge $e$ of
  $Q$ (the case $v \in \interior{Q}$ is treated analogously).  From
  Lemma~\ref{le:positive-eigenvalues} we know that there exist generators
  $p_i$, $p_j$ determining $v$ on opposite sides of $l$, the orthogonal line
  to the edge $e$ passing through $v$. From eq.~\eqref{eq:derivation} and
  Lemma~\ref{le:simple-facts-DC-SP}, we deduce that both $p_i$ and $p_j$ are
  centered. Now, for each generator $p_k$ with $v \in V_k$ in the same side
  of $l$ as $p_i$ (respectively $p_j$), we consider the triplet $(e,j,k)$
  (respectively $(e,i,k)$). Again resorting to eq.~\eqref{eq:derivation} and
  Lemma~\ref{le:simple-facts-DC-SP}, we conclude that $p_k$ is also centered.
  Finally, if a generator $p_k$ with $v \in V_k$ is such that $p_k \in l$,
  any of the triplets $(e,j,k)$ or $(e,i,k)$ can be invoked in a similar
  argument to ensure that $p_k$ is centered.  Therefore, $P \in A_{\DC} (Q)$,
  and hence $(Z_{-\LN (\partial \lg_{\VV}),\HHDC} \cap Q^n) \subset A_{\DC}
  (Q)$.  \quad
\end{proof}

\begin{proposition}\label{prop:distributed-gradient-SP}
  For the dynamical system~\eqref{eq:distributed-gradient-flow-SP}, the
  generators' location $P=(p_1,\dots,p_n)$ converges asymptotically to the
  largest weakly invariant set contained in the closure of $A_{\SP} (Q)=
  \setdef{ P \in Q^n}{i \in I(P) \implies p_i\in \IC (V_i)}$.
%
%   \sindex{ASP}{$A_{\SP}(Q)$}{Set of configurations $P \in Q^n$ where all
%     active generators are in the incenter set of its own Voronoi region}
%
\end{proposition}

\begin{proof}
  Let $a \in \setLieder{\LN (\partial \sm_{\VV})}{\HHSP}(P)$. By definition,
  $a= \LN [\sm_{\VV}] (P) \cdot \zeta$, for all $\zeta \in \partial \HHSP
  (P)$. Let $e \in \ED_{\SP}(\VV(P))$. If $e$ is an edge of type (a), i.e. a
  segment of the bisector determined by $p_i$ and $p_j$, we compute (cf.
  Proposition~\ref{prop:gradient-F}),
  \begin{align}\label{eq:derivation-SP}
    a & = \LN (\partial \sm_{\VV}) (P) \cdot \partial_{e} F_{i} \nonumber \\
    &= \LN (\partial \sm_{V_i (P)}) (P) \cdot \pi_i (\partial_{e} F_{i} ) +
    \LN ( \partial \sm_{V_j (P)}) (P) \cdot \pi_j (\partial_{e} F_{i} ) \, .
  \end{align}
  From Lemma~\ref{le:simple-facts-DC-SP}(iii) we conclude that $a \ge 0$, and
  the inequality is strict if either $p_i \not \in \IC (V_i)$ or $p_j \not
  \in \IC (V_j)$.  The same conclusion can be derived if $e$ is a vertex of
  type (b).  Therefore, $\max \setLieder{\LN (\partial \sm_{\VV})}{\HHSP}(P)
  \ge 0$ or $\setLieder{\LN (\partial \sm_{\VV})}{\HHSP}(P) = \emptyset$.
  Now, resorting to the LaSalle principle (Theorem~\ref{th:LaSalle}), we
  deduce that the solution $P:[0,+\infty) \rightarrow Q^n$ starting from
  $P_0$ converges to the largest weakly invariant set contained in
  $\ov{Z}_{\LN (\partial \sm_{\VV}),\HHSP} \cap \HHSP^{-1} (\le
  \HHSP(P_0),P_0)\cap Q^n$.  From eq.~\eqref{eq:derivation-SP}, and resorting
  to Proposition~\ref{prop:proj-generalized-gradients} and
  Lemma~\ref{le:simple-facts-DC-SP}, one can also show that $Z_{\LN (\partial
    \sm_{\VV}),\HHSP}\cap Q^n$ is equal to $A_{\SP} (Q)$.  \quad
\end{proof}

\begin{remark}
  {\rm The sets $A_{\DC}(Q)$ and $A_{\SP}(Q)$ are not closed in general. If
    $\dim Q = 1$, then it can be seen that they indeed are. In higher
    dimensions one can find sequences $\setdef{P_k\in Q^n}{k\in\natural}$
    in these sets which converge to configurations $P$ where not all active
    generators are centered.}
\end{remark}

% \begin{theorem}
%   The largest weakly invariant set under the
%   dynamics~\eqref{eq:distributed-gradient-flow-SP} contained in the
%   closure of $A_{\SP} (Q)$ is the set of all incenter Voronoi
%   configurations.
% \end{theorem}

\subsection{Distributed dynamical systems based on geometric centering}

Here, we propose alternative distributed dynamical systems for the
multi-center functions. Our design is directly inspired by the results in
Theorems~\ref{th:minima-DC} and~\ref{th:maxima-SP} on the critical points of
the multi-center functions $\HHDC$ and $\HHSP$.  For $i \in \until{n}$,
consider the dynamical systems
\begin{align}
  \dot{p}_i &= \CC (V_i) - p_i \, , \label{eq:circumcenter-flow} \\
  \dot{p}_i & \in \IC (V_i) - p_i \, . \label{eq:incenter-flow}
\end{align}
Alternatively, we may write $\dot{P} = \CC (\VV (P)) - P$ and $\dot{P} \in
\IC (\VV (P)) - P$. Note that both systems are Voronoi-distributed. Also,
note that the vector field~\eqref{eq:circumcenter-flow} is continuous, since
the circumcenter of a polygon depends continuously on the location of its
vertexes, and the location of the vertexes of the Voronoi partition depends
continuously on the location of the generators; see~\cite{AO-BB-KS-SNC:00}.
However, eq.~\eqref{eq:incenter-flow} is a differential inclusion, since the
incenter sets may not be singletons.  By Lemma~\ref{le:existence-solution},
the existence of solutions to eq.~\eqref{eq:incenter-flow} is guaranteed by
the following result.

\begin{proposition}\label{le:move-toward-incenter-well-defined}
  Consider the set-valued map $\IC(\VV) - \Id: Q^n \rightarrow
  2^{(\real^2)^n}$ given by $P \mapsto \IC (\VV (P)) - P$. Then $\IC(\VV) -
  \Id$ is upper semicontinuous with nonempty, compact and convex values.
\end{proposition}

\begin{proof}
  Clearly, the map $\IC(\VV) - \Id$ takes nonempty and compact values. From
  Lemma~\ref{le:incenter-set-convex-segment}, we also know that it takes
  convex values. Furthermore, since the identity map is continuous, it
  suffices to check that $P \mapsto \IC(\VV (P))$ is upper semicontinuous.
  We then have to verify that, given $P_0 \in Q^n$, for each $\eps>0$,
  there exists $\delta >0$ such that
  \begin{align}\label{eq:checking-upper-semi-continuous}
    \IC (\VV (P)) \subset \IC (\VV (P_0)) + B_{2n} (0,\eps) \, , \quad
    \text{if} \; \, \| P - P_0 \| \le \delta \, .
  \end{align} 
  Now, for each $i$, if $\IC(V_i(P_0))$ is not a singleton, then it is a
  segment (cf. Lemma~\ref{le:incenter-set-convex-segment}) whose extremal
  points $q_{i1}(P_0)$, $q_{i2}(P_0)$ are the intersection points of some
  bisectors of the edges of the Voronoi cell. It is clear that
  $q_{i\alpha}(P) \rightarrow q_{i\alpha}(P_0)$ when $P \rightarrow P_0$
  for $\alpha=1,2$. Therefore, given $\eps>0$, one can choose $\delta_i >0
  $ such that if $\| P - P_0 \| \le \delta_i$, then $\|q_{i\alpha}(P) -
  q_{i\alpha}(P_0) \| \le \eps/n$.  Since $\IC (V_i(P))$ is contained in
  the segment joining $q_{i1} (P)$ and $q_{i2}(P)$, we deduce $\IC (V_i(P))
  \subset \IC (V_i (P_0)) + B_2 (0,\eps/n)$.  On the other hand, if
  $\IC(V_i(P_0))$ is a singleton, then it coincides with the intersection
  points $q_{i1}(P_0), \dots, q_{im}(P_0)$ of some bisectors of the edges
  of the Voronoi cell. The above reasoning also guarantees that there exits
  $\delta_i>0$ such that $q_{i\alpha} (P) \in \IC (V_i (P_0)) + B_2
  (0,\eps/n)$, $\alpha=1,\dots,m$, if $\| P - P_0 \| \le \delta_i$. Since
  $\IC (V_i(P))$ is contained in one of the segments joining the points
  $q_{i1} (P), \dots, q_{im}(P)$, we again deduce $\IC (V_i(P)) \subset \IC
  (V_i (P_0)) + B_2 (0,\eps/n)$.  The statement
  in~\eqref{eq:checking-upper-semi-continuous} follows by taking the
  minimum of $\delta_1, \dots, \delta_n$.  \quad
\end{proof}

Having established the existence of solutions, one can also see that the
compact set $Q^n$ is strongly invariant for the vector field $\CC(\VV) -
\Id$ and for the differential inclusion $\IC(\VV) - \Id$.  Next, we
characterize the asymptotic convergence of the dynamical systems under
study.

\begin{proposition}\label{prop:move-to-center}
  For the dynamical system~\eqref{eq:circumcenter-flow}
  (respectively~\eqref{eq:incenter-flow}), the generators' location
  $P=(p_1,\dots,p_n)$ converges asymptotically to the largest weakly
  invariant set contained in the closure of $A_{\DC} (Q)$ (respectively  in the
  closure of $A_{\SP} (Q)$).
\end{proposition}

\begin{proof}
  The proof of this result is parallel to the proofs of
  Propositions~\ref{prop:distributed-gradient-DC}
  and~\ref{prop:distributed-gradient-SP}. The sequence of steps is the same
  as before, though now one resorts to Lemma~\ref{le:simple-facts-DC-SP}(ii)
  and Lemma~\ref{le:simple-facts-DC-SP}(iv). The only additional observation
  is that, when computing the set-valued Lie derivative for
  eq.~\eqref{eq:incenter-flow}, one has that $a \in \setLieder{\IC(\VV) -
    \Id}{\HHSP}(P)$ if and only if there exists $x \in \IC(\VV(P))$ such that
  $a = (x-P) \cdot \zeta$, for any $\zeta \in \partial \HHSP (P)$. The
  application of Lemma~\ref{le:simple-facts-DC-SP} guarantees that $a \ge 0$,
  and that the inequality is strict if any of the active generators is not in
  its corresponding incenter set. \quad
\end{proof}

% \begin{theorem}
%   The largest weakly invariant set under the
%   dynamics~\eqref{eq:circumcenter-flow} contained in the closure of
%   $A_{\DC} (Q)$ is the set of all circumcenter Voronoi configurations.
% \end{theorem}
% \begin{theorem}
%   The largest weakly invariant set under the
%   dynamics~\eqref{eq:incenter-flow} contained in the closure of
%   $A_{\SP} (Q)$ is the set of all incenter Voronoi configurations.
% \end{theorem}

\subsection{Simulations}

\begin{figure*}[htbp]
  \centering%  
  \fbox{\includegraphics[width=.3\linewidth]{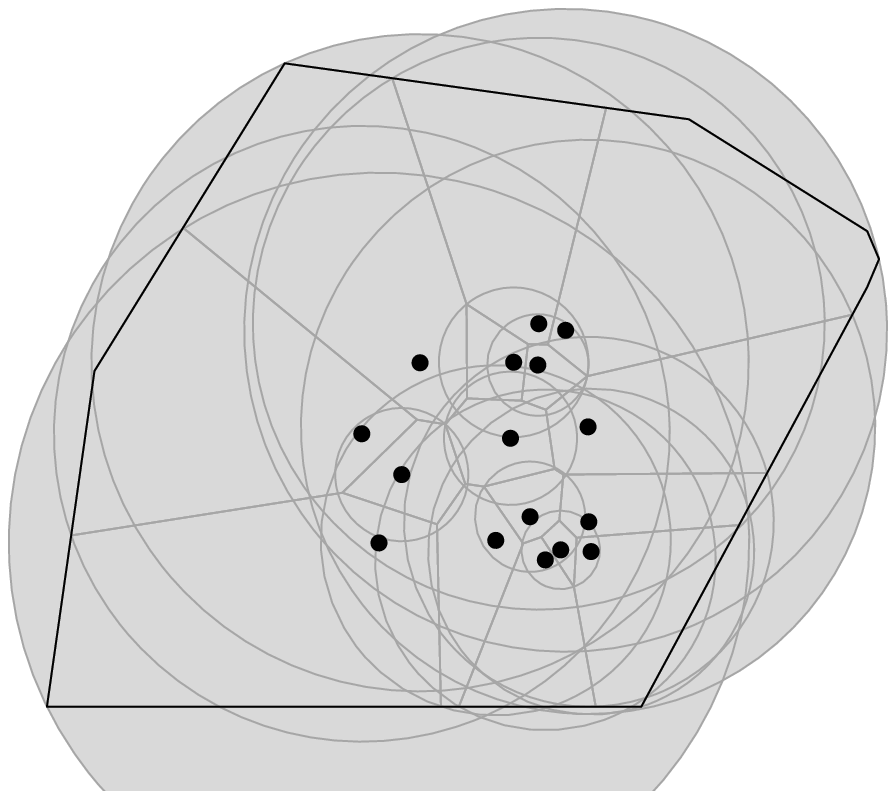}}
  \fbox{\includegraphics[width=.3\linewidth]{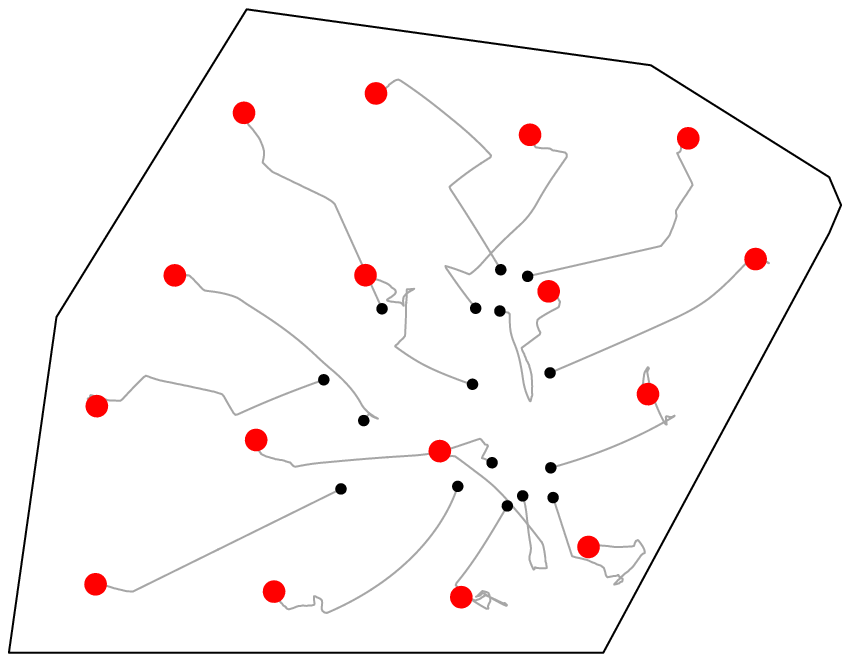}}
  \fbox{\includegraphics[width=.3\linewidth]{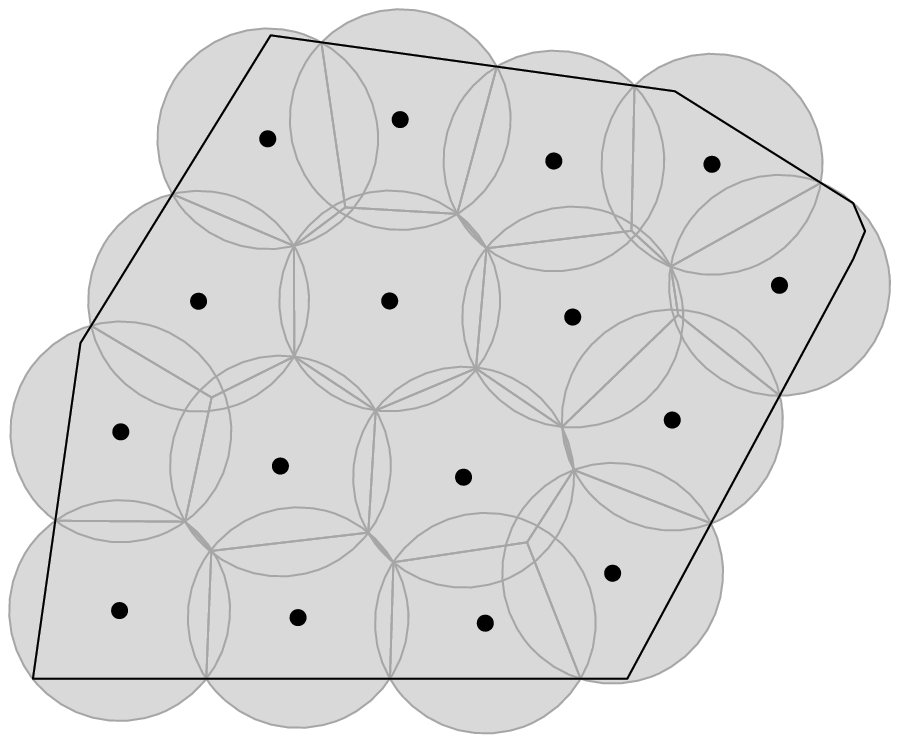}}
  \caption{``Toward the furthest'' algorithm for $16$ generators
    in a convex polygonal environment.  The left (respectively, right)
    figure illustrates the initial (respectively, final) locations and
    Voronoi partition. The central figure illustrates the network
    evolution. After $2$ seconds, the multi-center function is
    approximately $.39504$ meters.}
  \label{fig:coverage-toward-furthest}
\end{figure*}

\begin{figure*}[htbp]
  \centering
  \fbox{\includegraphics[width=.3\linewidth]{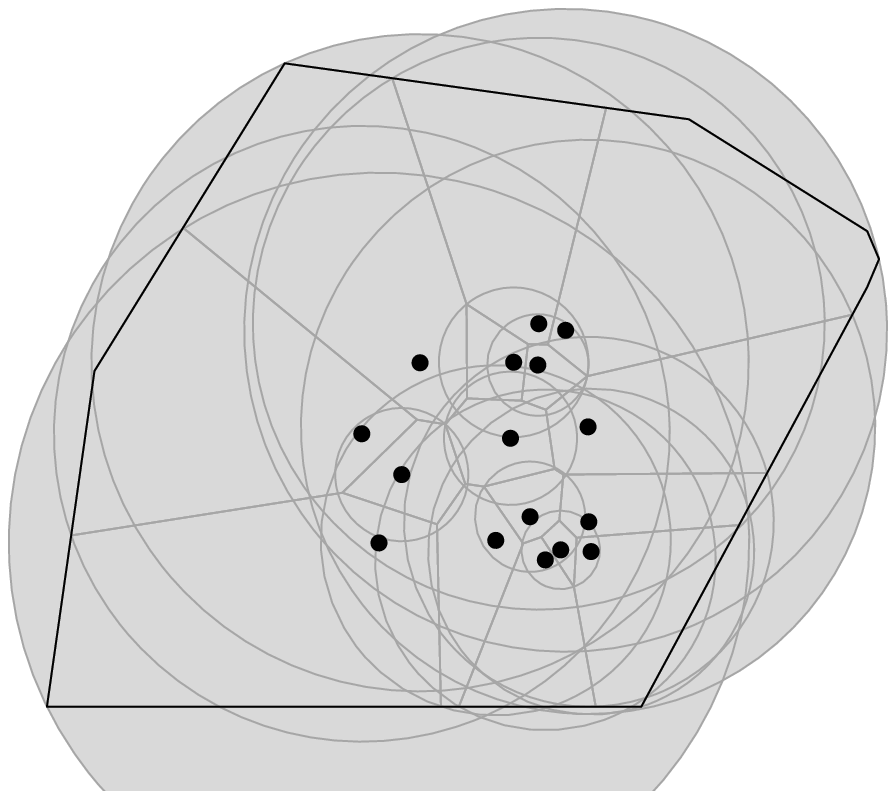}}
  \fbox{\includegraphics[width=.3\linewidth]{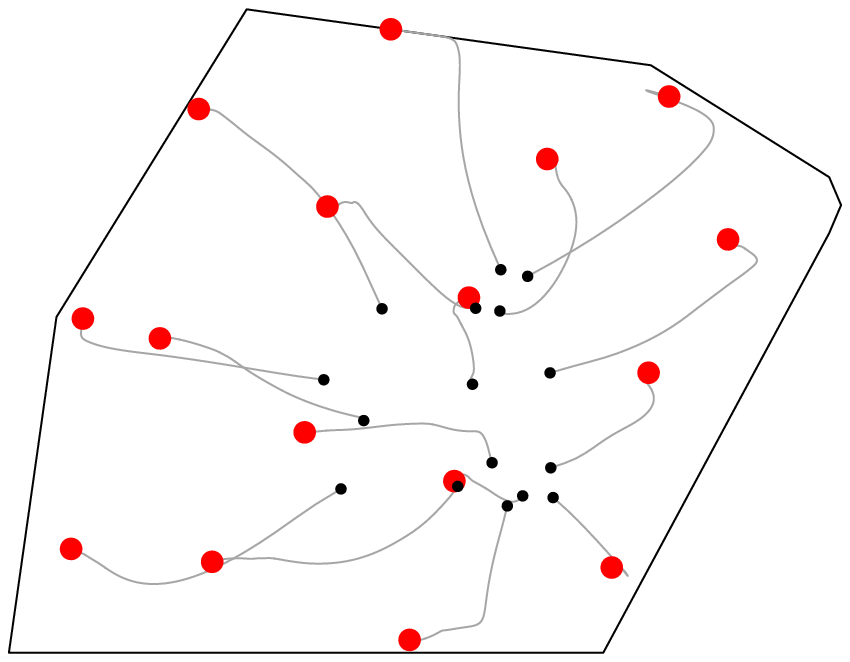}}
  \fbox{\includegraphics[width=.3\linewidth]{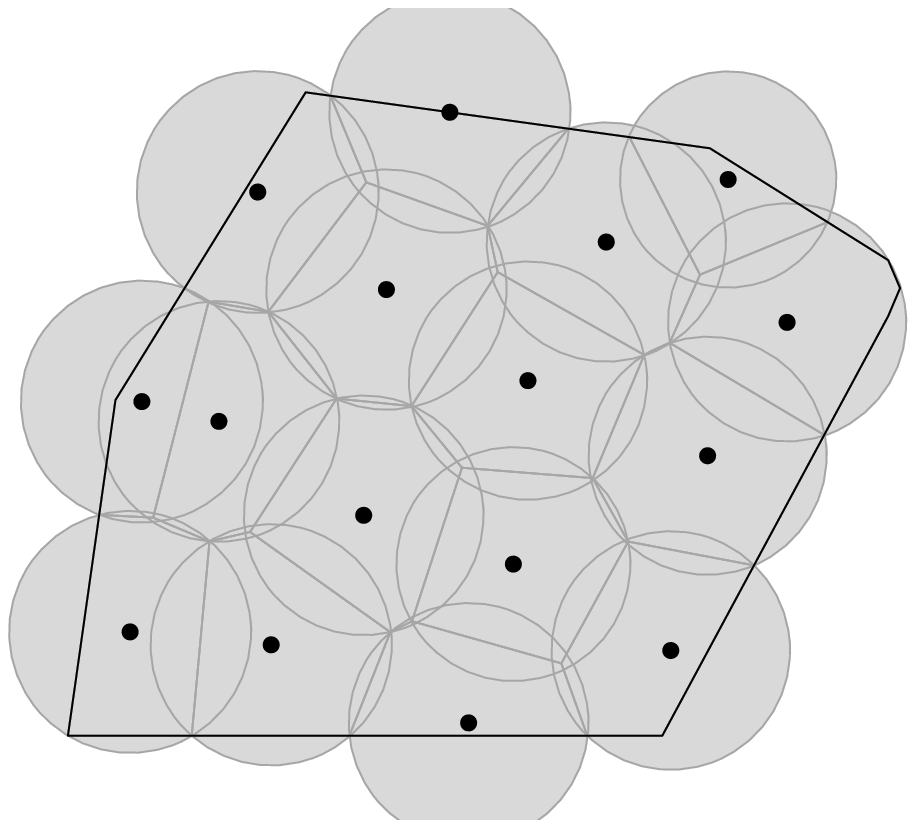}}
  \caption{``Move-toward-the-circumcenter'' algorithm for $16$ generators
    in a convex polygonal environment.  The left (respectively, right)
    figure illustrates the initial (respectively, final) locations and
    Voronoi partition. The central figure illustrates the network
    evolution.  After $20$ seconds, the multi-center function is
    approximately $0.43273$ meters.}
  \label{fig:coverage-lloyd-DC}
\end{figure*}

\begin{figure*}[htbp]
  \centering
  \fbox{\includegraphics[width=.3\linewidth]{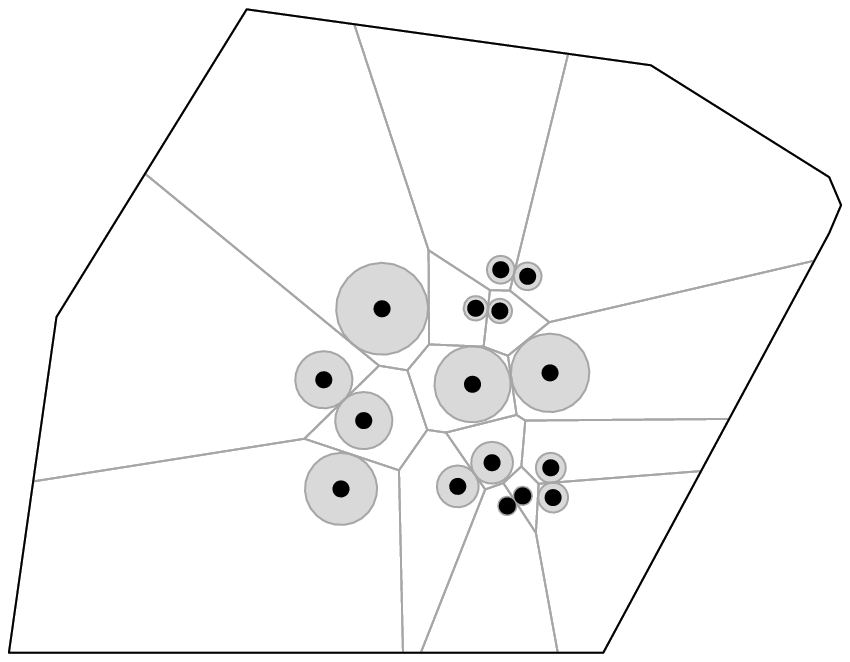}}
  \fbox{\includegraphics[width=.3\linewidth]{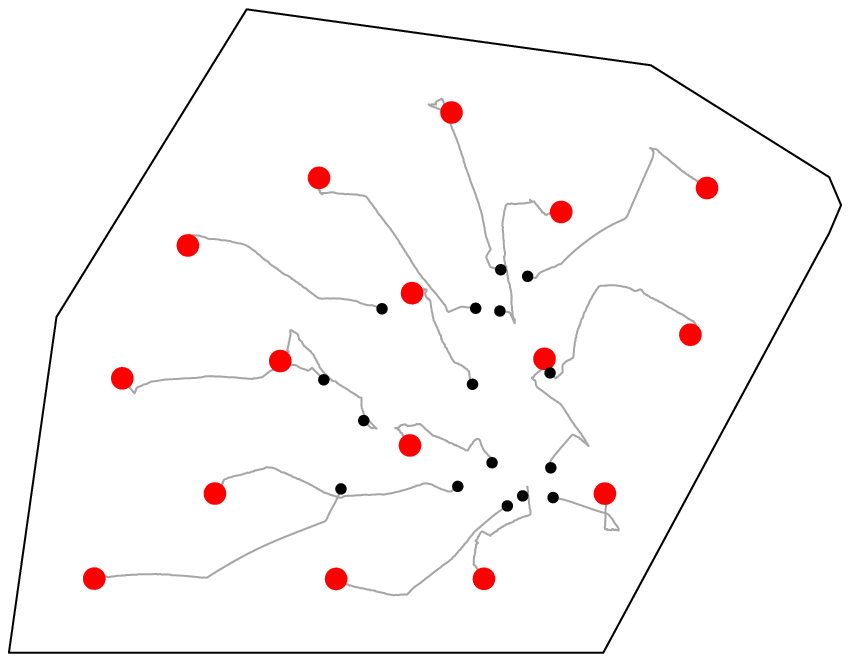}}
  \fbox{\includegraphics[width=.3\linewidth]{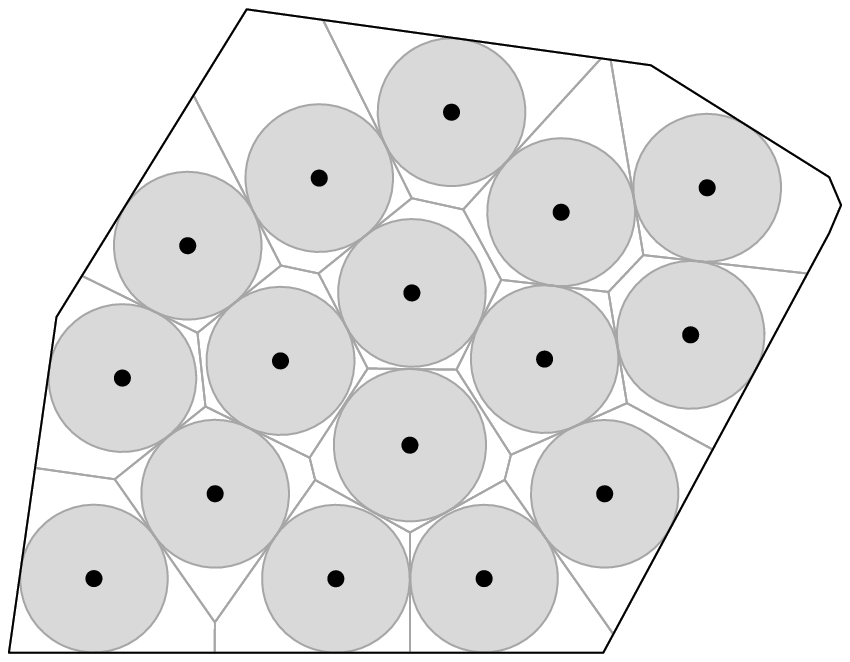}}
  \caption{``Away-from-closest'' algorithm for $16$  generators
    in a convex polygonal environment.  The left (respectively, right)
    figure illustrates the initial (respectively, final) locations and
    Voronoi partition. The central figure illustrates the network
    evolution.  After $2$ seconds, the multi-center function is
    approximately $.26347$ meters.}
  \label{fig:away-from-closest}
\end{figure*}

\begin{figure*}[htbp]
  \centering
  \fbox{\includegraphics[width=.3\linewidth]{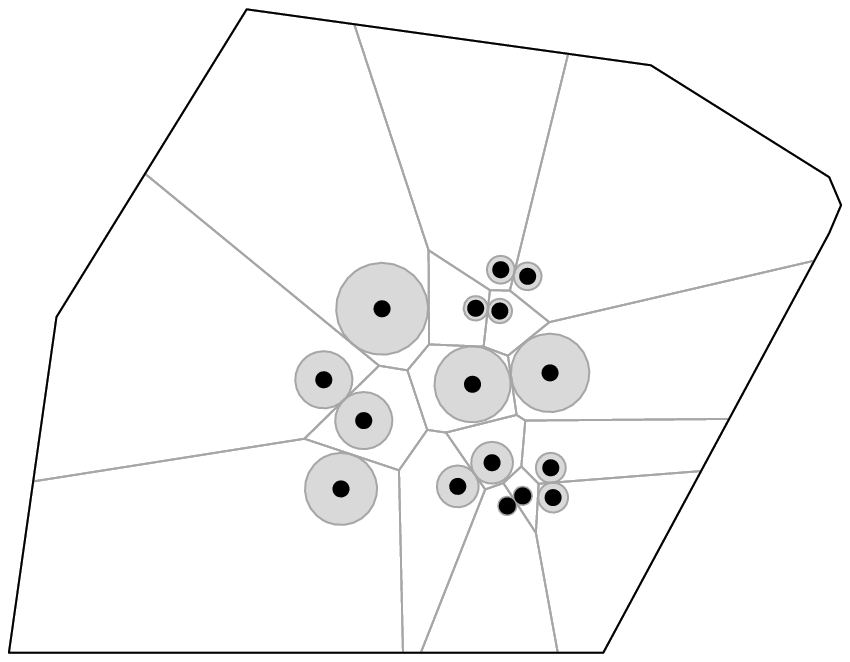}}
  \fbox{\includegraphics[width=.3\linewidth]{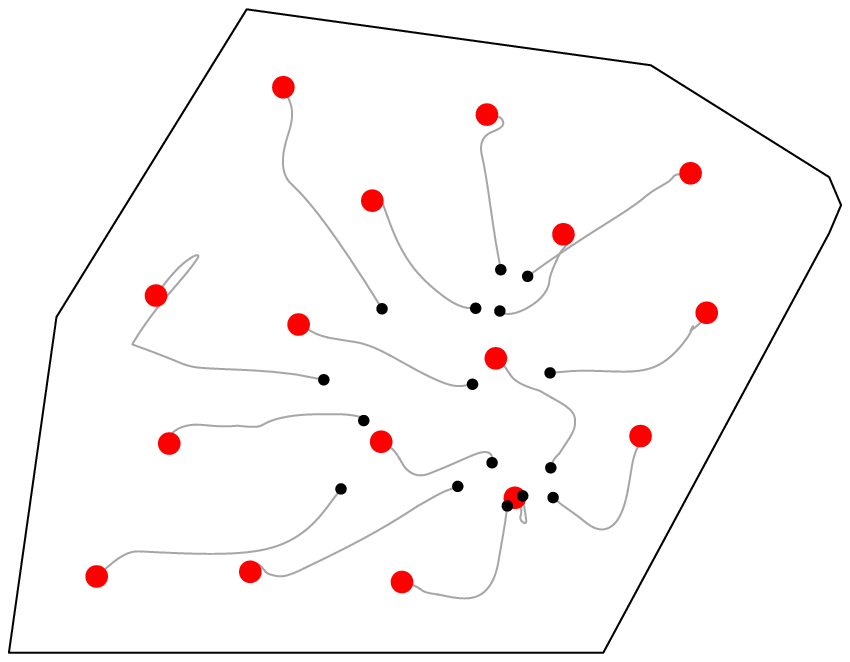}}
  \fbox{\includegraphics[width=.3\linewidth]{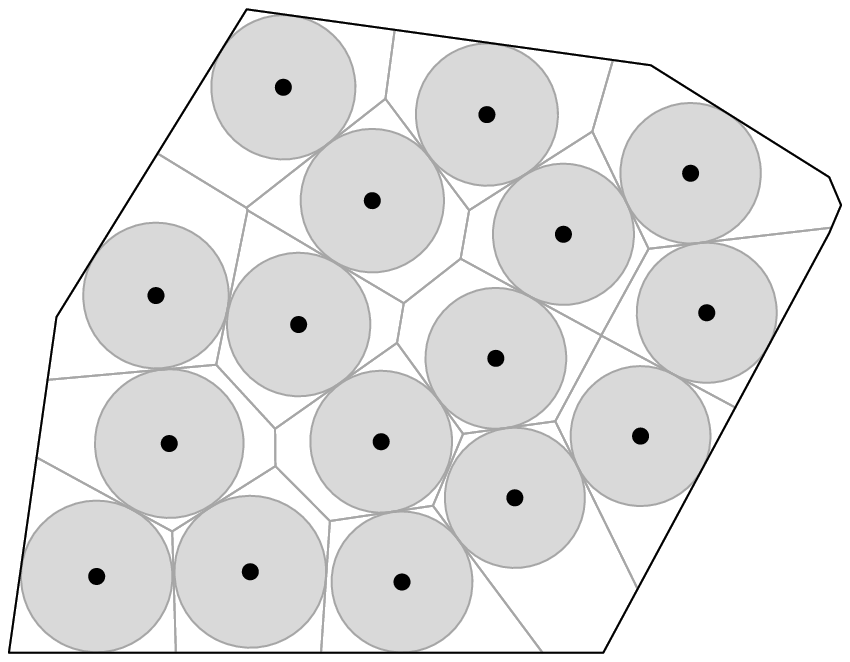}}
  \caption{``Move-toward-the-incenter'' algorithm for $16$ generators  in a convex
    polygonal environment.  The left (respectively, right) figure
    illustrates the initial (respectively, final) locations and Voronoi
    partition. The central figure illustrates the network evolution.  After
    $20$ seconds, the multi-center function is approximately $.2498$
    meters.}
  \label{fig:coverage-lloyd-SP}
\end{figure*}

To illustrate the performance of the distributed coordination algorithms,
we include some simulation results. The algorithms are implemented in
\texttt{Mathematica} as a single centralized program.  We compute the
bounded Voronoi diagram of a collection of points using the
\texttt{Mathematica} package \texttt{ComputationalGeometry}. We compute the
circumcenter of a polygon via the algorithm in~\cite{SS:91} and the
incenter set via the \texttt{LinearProgramming} solver in
\texttt{Mathematica}.  Measuring displacements in meters, we consider the
domain determined by the vertexes
\begin{align*}
  \{(0,0), (2.5,0), (3.45,1.5), (3.5,1.6), (3.45,1.7), (2.7,2.1), (1.,2.4),
  (.2,1.2)\}.
\end{align*}
In Figs.~\ref{fig:coverage-toward-furthest}
and~\ref{fig:coverage-lloyd-DC}, we illustrate the performance of the
dynamical systems~\eqref{eq:distributed-gradient-flow-DC}
and~\eqref{eq:circumcenter-flow}, respectively, minimizing the
multi-circumcenter function $\HHDC$.  In Figs.~\ref{fig:away-from-closest}
and~\ref{fig:coverage-lloyd-SP}, we illustrate the performance of the
dynamical systems~\eqref{eq:distributed-gradient-flow-SP}
and~\eqref{eq:incenter-flow}, respectively, maximizing the multi-incenter
function $\HHSP$.  Observing the final configurations in the four figures,
one can verify, visually and numerically, that the active generators are
asymptotically centered as forecast by our analysis.

\section{Conclusions}\label{sec:conclusions}

We have introduced two multi-center functions that provide
quality-of-service measures for mobile networks. We have shown that both
functions are globally Lipschitz and regular, and we have computed their
generalized gradients. Furthermore, under certain technical conditions, we
have characterized via nonsmooth analysis their critical points as center
Voronoi configurations and as solutions of disk-covering and sphere-packing
problems. We have also considered various algorithms that extremize the
multi-center functions.  First, we considered the nonsmooth gradient flows
induced by their respective generalized gradients.  Second, we devised a
novel strategy based on the generalized gradients of the 1-center functions
of each generator. Third, we introduced and characterized a geometric
centering strategy with resemblances to the classical Lloyd algorithm.  We
have unveiled the remarkable geometric interpretations of these algorithms,
discussed their distributed character and analyzed their asymptotic
behavior using nonsmooth stability analysis.

Future directions of research include: (i) sharpening the asymptotic
convergence results for the proposed dynamical systems, (ii)
considering the setting of convex polytopes in $\real^N$, for $N>2$,
(iii) understanding in what sense the proposed multi-circumcenter and
the multi-incenter problems can be shown to be dual, and (iv)
analyzing other meaningful geometric optimization problems and their
relations with cooperative behaviors.

\clearpage

\begin{symbollist}
\item[\textbf{Symbol}] \textbf{Description and page(s) when
    applicable} \vspace{3\lineskip}
\item[$A_{\DC}(Q)$] Set of configurations $P \in Q^n$ where all active
  generators are in the circumcenter of its own Voronoi region, 25
\item[$A_{\SP}(Q)$] Set of configurations $P \in Q^n$ where all active
  generators are in the incenter set of its own Voronoi region, 26
\item[$\CC(Q)$] Circumcenter of polytope $Q$, 6
\item[$\CR(Q)$] Circumradius of polytope $Q$, 6
\item[$\D_S$] Distance function to the convex set $S$, 4
\item[$\ED(Q)$] Edges of polygon $Q$, 4
\item[$\ED_{\SP}(\VV(P))$] Edges where the value of $\HHSP(P)$ is   attained, 6
\item[$e(i)$] Edge of $\VV (P)$ belonging to $V_i$ and to the boundary of   $Q$, 5
\item[$e(i,j)$] Edge of $\VV (P)$ determined by $p_i$ and $p_j$, 5
\item[$F_i(P)$] Smallest distance from $p_i$ to the boundary of $V_i (P)$, 14
\item[$G_i(P)$] Largest distance from $p_i$ to the boundary of $V_i (P)$, 14
\item[$\partial f$] Generalized gradient of the   locally Lipschitz function $f$, 7
\item[$\HHDC$] Multi-circumcenter function, 6
\item[$\HHSP$] Multi-incenter function, 6
\item[$\IC(Q)$] Incenter set of polytope $Q$, 6
\item[$\IR(Q)$] Inradius of polytope $Q$, 6
\item[{$K[X]$}] Filippov mapping associated with a measurable and
  essentially locally bounded mapping $X:\real^N \rightarrow \real^N$
  , 8
\item[$\lambda(e,i,j)$] Scalar function associated with the vertex
  $v(e,i,j)$, 16
\item[{$\LN (S)$}] Least-norm element of the convex set $S$, 7
\item[$\lg_Q(p)$] Largest distance from $p$ to the boundary of $Q$, 10
\item[$\mu(i,j,k)$] Scalar function associated with the vertex   $v(i,j,k)$, 17
\item[$\NN (P,i)$, $\NN(i)$] Set of neighbors of the $i$th generator
  at configuration $P$, 4
\item[$n_{e(i,j)}$] Unit normal to $e(i,j)$  pointing toward  $\interior{V_i(P)}$, 5
\item[$n_{e(i)}$] Unit normal to $e(i)$ pointing toward
  $\interior{Q}$, 5
\item[$\proj{S}$] Orthogonal projection onto the convex set $S$, 4
\item[$\pi_i$] Canonical projection from $Q^n$ onto the $i$th factor,
  4
\item[$\setLieder{X}{f}$] Set-valued Lie derivative of $f$ with respect to $X$, 8
\item[$\sm_Q(p)$] Smallest distance from $p$ to the boundary of $Q$,
  10
\item[$v(i,j,k)$] Vertex of $\VV (P)$ determined by $p_i$, $p_j$ and   $p_k$, 4
\item[$v(e,i,j)$] Vertex of $\VV (P)$ determined by $e \in \ED (Q)$
  and $p_i$, $p_j$, 4
\item[$v(e,f,i)$] Vertex of $\VV (P)$ determined by $e,f \in \ED (Q)$   and $p_i$, 4
\item[$\VE_{\DC}(\VV(P))$] Vertexes of $\VV (P)$ where the value of   $\HHDC (P)$ is attained, 6
\item[$\versus{v}$] Unit vector in the direction of $0 \neq v \in
  \real^N$, 4
\item[$\VE(Q)$] Vertexes of polygon $Q$, 4
\item[$\VV(P)$] Voronoi partition of $Q$ generated by   $P=(p_1,\dots,p_n)$, 4
\item[$Z_{X,f}$] Set formed by points $x \in\real^N$ such that $0$
  belongs to $\setLieder{X}{f}(x)$, 9

\end{symbollist}


\begin{thebibliography}{10}

\bibitem{PKA-MS:98}
{\sc P.~K. Agarwal and M.~Sharir}, {\em Efficient algorithms for geometric
  optimization}, ACM Computing Surveys, 30 (1998), pp.~412--458.

\bibitem{RCA:98}
{\sc R.~C. Arkin}, {\em Behavior-Based Robotics}, Cambridge University Press,
  New York, NY, 1998.

\bibitem{AB-FC:99}
{\sc A.~Bacciotti and F.~Ceragioli}, {\em Stability and stabilization of
  discontinuous systems and nonsmooth {L}yapunov functions}, {ESAIM.} Control,
  Optimisation \& Calculus of Variations, 4 (1999), pp.~361--376.

\bibitem{DPB-JNT:97}
{\sc D.~P. Bertsekas and J.~N. Tsitsiklis}, {\em Parallel and Distributed
  Computation: Numerical Methods}, Athena Scientific, 1997.

\bibitem{VB-HM-VS:99}
{\sc V.~Boltyanski, H.~Martini, and V.~Soltan}, {\em Geometric methods and
  optimization problems}, vol.~4 of Combinatorial optimization, Kluwer Academic
  Publishers, Dordrecht, Boston, 1999.

\bibitem{SB-LV:02}
{\sc S.~Boyd and L.~Vandenberghe}, {\em Convex optimization}.
\newblock Preprint, Dec. 2002.

\bibitem{RWB:91}
{\sc R.~W. Brockett}, {\em Dynamical systems that sort lists, diagonalize
  matrices, and solve linear programming problems}, Linear Algebra and its
  Applications, 146 (1991), pp.~79--91.

\bibitem{HC:99}
{\sc H.~Choset}, {\em Nonsmooth analysis, convex analysis, and their
  applications to motion planning}, International Journal of Computational
  Geometry and Applications, 9 (1999), pp.~447--469.

\bibitem{FHC:83}
{\sc F.~H. Clarke}, {\em Optimization and Nonsmooth Analysis}, Canadian
  Mathematical Society Series of Monographs and Advanced Texts, John Wiley \&
  Sons, 1983.

\bibitem{JC-SM-TK-FB:02j}
{\sc J.~Cort{\'e}s, S.~Mart{\'\i}nez, T.~Karatas, and F.~Bullo}, {\em Coverage
  control for mobile sensing networks}, IEEE Transactions on Robotics and
  Automation,  (2002).
\newblock Conditionally accepted.

\bibitem{MdB-MvK-MO:97}
{\sc M.~de~Berg, M.~van Kreveld, and M.~Overmars}, {\em Computational Geometry:
  Algorithms and Applications}, Springer Verlag, New York, NY, 1997.

\bibitem{JPD-JPO-VK:01}
{\sc J.~P. Desai, J.~P. Ostrowski, and V.~Kumar}, {\em Modeling and control of
  formations of nonholonomic mobile robots}, IEEE Transactions on Robotics and
  Automation, 17 (2001), pp.~905--8.

\bibitem{ZD:95}
{\sc Z.~Drezner}, ed., {\em Facility Location: A Survey of Applications and
  Methods}, Springer Series in Operations Research, Springer Verlag, New York,
  NY, 1995.

\bibitem{QD-VF-MG:99}
{\sc Q.~Du, V.~Faber, and M.~Gunzburger}, {\em Centroidal {V}oronoi
  tessellations: applications and algorithms}, SIAM Review, 41 (1999),
  pp.~637--676.

\bibitem{AFF:88}
{\sc A.~F. Filippov}, {\em Differential Equations with Discontinuous Righthand
  Sides}, vol.~18 of Mathematics and Its Applications, Kluwer Academic
  Publishers, Dordrecht, 1988.
\newblock Original Russian edition: \textit{Differentsial'nye Uravneniya s
  Razryvnoi Pravoi Chast'yu}, Nauka, Moscow, 1985.

\bibitem{RMG-DLN:98}
{\sc R.~M. Gray and D.~L. Neuhoff}, {\em Quantization}, IEEE Transactions on
  Information Theory, 44 (1998), pp.~2325--2383.
\newblock Commemorative Issue 1948-1998.

\bibitem{UH-JBM:94}
{\sc U.~Helmke and J.~Moore}, {\em Optimization and Dynamical Systems},
  Springer Verlag, New York, NY, 1994.

\bibitem{AJ-JL-ASM:02}
{\sc A.~Jadbabaie, J.~Lin, and A.~S. Morse}, {\em Coordination of groups of
  mobile autonomous agents using nearest neighbor rules}, IEEE Transactions on
  Automatic Control,  (2003).
\newblock To appear.

\bibitem{NEL-EF:01}
{\sc N.~E. Leonard and E.~Fiorelli}, {\em Virtual leaders, artificial
  potentials, and coordinated control of groups}, in {IEEE} Conf. on Decision
  and Control, Orlando, FL, Dec. 2001, pp.~2968--2973.

\bibitem{YL-KMP-MMP:00}
{\sc Y.~Liu, K.~M. Passino, and M.~M. Polycarpou}, {\em Stability analysis of
  m-dimensional asynchronous swarms with a fixed communication topology}, IEEE
  Transactions on Automatic Control, 48 (2003), pp.~76--95.

\bibitem{SPL:82}
{\sc S.~P. Lloyd}, {\em Least squares quantization in {PCM}}, IEEE Transactions
  on Information Theory, 28 (1982), pp.~129--137.
\newblock Presented as Bell Laboratory Technical Memorandum at a 1957 Institute
  for Mathematical Statistics meeting.

\bibitem{DGL:84}
{\sc D.~G. Luenberger}, {\em Linear and Nonlinear Programming}, Addison-Wesley,
  Reading, Massachusetts, second~ed., 1984.

\bibitem{AO-BB-KS-SNC:00}
{\sc A.~Okabe, B.~Boots, K.~Sugihara, and S.~N. Chiu}, {\em Spatial
  Tessellations: Concepts and Applications of Voronoi Diagrams}, Wiley Series
  in Probability and Statistics, John Wiley \& Sons, New York, NY, second~ed.,
  2000.

\bibitem{ROS-RMM:03b}
{\sc R.~Olfati-Saber and R.~M. Murray}, {\em Agreement problems in networks
  with directed graphs and switching topology}, in {IEEE} Conf. on Decision and
  Control, 2003.
\newblock Submitted.

\bibitem{BP-SSS:87}
{\sc B.~Paden and S.~S. Sastry}, {\em A calculus for computing {F}ilippov's
  differential inclusion with application to the variable structure control of
  robot manipulators}, IEEE Transactions on Circuits and Systems, 34 (1987),
  pp.~73--82.

\bibitem{JMR-GTT:90}
{\sc J.-M. Robert and G.~T. Toussaint}, {\em Computational geometry and
  facility location}, in Proc. International Conf. on Operations Research and
  Management Science, vol.~B, Manila, The Philippines, Dec. 1990, pp.~1--19.

\bibitem{DS-BP:94}
{\sc D.~Shevitz and B.~Paden}, {\em Lyapunov stability theory of nonsmooth
  systems}, IEEE Transactions on Automatic Control, 39 (1994), pp.~1910--1914.

\bibitem{SS:91}
{\sc S.~Skyum}, {\em A simple algorithm for computing the smallest circle},
  Information Processing Letters, 37 (1991), pp.~121--125.

\bibitem{AS-ZD:96}
{\sc A.~Suzuki and Z.~Drezner}, {\em The $p$-center location problem in an
  area}, Location Science, 4 (1996), pp.~69--82.

\bibitem{IS-MY:99}
{\sc I.~Suzuki and M.~Yamashita}, {\em Distributed anonymous mobile robots:
  Formation of geometric patterns}, SIAM Journal on Computing, 28 (1999),
  pp.~1347--1363.

\bibitem{HT-AJ-GJP:03b}
{\sc H.~Tanner, A.~Jadbabaie, and G.~J. Pappas}, {\em Stable flocking of mobile
  agents, {Part II:} dynamic topology}, in {IEEE} Conf. on Decision and
  Control, Maui, Hawaii, Dec. 2003.
\newblock Submitted.

\end{thebibliography}
\end{document}